\documentclass[a4paper,12pt]{article}

\def\displayandname#1{\rlap{$\displaystyle\csname #1\endcsname$}%
                      \qquad \texttt{\char92 #1}}

\newtheorem{thm}{Theorem}[section]
\newtheorem{pro}[thm]{Proposition}
\newtheorem{lem}[thm]{Lemma}

\newtheorem{cor}[thm]{Corollary}
\newtheorem{fac}[thm]{Fact}

\newtheorem{exa}[thm]{Example}
\newtheorem{df}[thm]{Definition}
\newtheorem{rem}[thm]{Remark}
\newtheorem{nt}[thm]{Note}

\begin{document}
\bibliographystyle{acm}
\title{Statistics on Wreath Products}
\author{Michael Fire \thanks{This paper is a part of the author's Master Thesis, which was
written at Bar-Ilan University under the supervision of R. M. Adin
and Y. Roichman.}\
\thanks{Partially supported by EC's IHRP Programme, within the
Research Training Network ``Algebraic Combinatorics
in Europe'', grant HPRN-CT-2001-00272}\\
Department of Mathematics\\
Bar Ilan University\\
Ramat Gan 52900, Israel\\
{\em mickyfi@gmail.com}}
\maketitle

\section{Introduction}
    The \emph{colored permutation groups} are fundamental objects in
    much of today's mathematics. A better understanding of these
    groups may help to advance research in many fields. One method of
    studying these groups is by using numerical statistics and finding their
    generating functions. This method was successfully applied in the
    case of "one-colored" permutations groups, the \emph{symmetric groups}.
    MacMahon \cite{MM} considered four different statistics for a
    permutation $\pi$ in the \emph{symmetric group}: the number of \emph{descents}
    ($des(\pi)$), the number of \emph{excedance} ($exc(\pi)$), the
    \emph{length} statistic ($\ell(\pi)$), and the \emph{major index}
    ($maj(\pi)$). MacMahon showed, algebraically, that \emph{excedance}
    number is equidistributed with \emph{descent} number, and that
    \emph{length} is equidistributed with \emph{major index} over the \emph{symmetric
    groups}.\\

    When we talk about permutation statistics we generally discuss
    about two main types of statistics: \emph{Eulerian} statistics,
    which are equidistributed with the \emph{descent} number,
    and \emph{Mahonian} statistics which are equidistributed with \emph{length}.
    Through the years many generalizations to MacMahon's results were
    found \cite{Den,DF,FS,FZ} using combinations of statistics which
    are equidistributed with each other over the \emph{symmetric
    groups}.\\

    Recently, Adin and Roichman \cite{AR} generalized MacMahon's result
    on the \emph{major index} to "two-colored" permutation groups, the \emph{signed permutations groups}, by
    introducing a new Mahonian statistic, the
    \emph{flag major index}. Bagno \cite{BE} introduced
    a new Mahonian statistic, \emph{lmaj}, which is equidistributed with
    \emph{length} over general \emph{colored permutation groups}. These results
    open the gate for trying to generalize results that were obtained
    on \emph{symmetric groups} to general \emph{colored
    permutations groups}. In this paper we attempt to go in the above
    mentioned direction, and generalize known theorems from the
    \emph{symmetric groups} to the \emph{signed permutations groups}, and to general
    \emph{colored permutation groups}.\\

    The paper is organized as follows: we start by giving necessary background in Section 2. In
    Section 3 we present our main results. In Section 4 we study permutation
    statistics with respect to different linear orders, we prove that the
    \emph{flag major index} is equidistributed with
    \emph{length} over the \emph{signed permutation
    groups} for every linear order, and also find a large collection of linear orders, which
    \emph{sign} and \emph{flag major index} equidistributed with each other.
    In Section 5 we present a new method for calculating generating
    functions on the \emph{signed permutations groups} in the
    \emph{natural order}:
    \[N: \  -n<-(n-1)<\ldots<-1<1<\ldots<n,\]
    and also calculate generating functions on a well-known
    subgroup of the \emph{signed permutation group}: $D_n$ (to be
    defined below). In Section 6 we move to calculating generating
    functions for statistics, but this time on general \emph{colored permutation
    groups} in the \emph{friends order}:
    \[F: \ 1^{[r-1]}<\ldots<1<2^{[r-1]}<\ldots<2<\ldots<n^{[r-1]}<\ldots<n.\]
    We present a new method of generalizing equidistributed statistics over the \emph{symmetric
    groups} to general \emph{colored permutations groups}. Using this method we find
    new Mahonian and Eulerian statistics, and generalize known
    theorems due to Foata; Zeilberger and Sch\"utzenberger. We conclude the paper in Section
    7, where we introduce the \emph{flag inversion number} and study its
    properties.

\section{Background}
    \subsection{Statistics on the Symmetric Group}
        In this subsection we present the main definitions, notation, and
        theorems on the \emph{symmetric groups} (i.e., the Weyl groups of type A),
        denoted $S_n$.

        \begin{df} Let be \textbf{N} the set of all the \emph{natural numbers}, a \emph{permutation} of order $n \in \textbf{N}$ is a bijection
            $\pi:\{1,2,3,\ldots,n\}\to \{1,2,3,\ldots,n\} $.
        \end{df}

        \begin{rem} {Permutations are traditionally written in a
        two-line notation of:
            \[\pi=
            \left(%
                \begin{array}{ccccc}
                    1 & 2 & 3 & \ldots & n \\
                    \pi(1)& \pi(2) & \pi(3) & \ldots & \pi(n) \\
                \end{array}
            \right),\] however for convenience we will use the shorter
            notation:
            \[\pi=[\pi(1),\pi(2),\pi(3),\ldots,\pi(n)].\]
            For example: $\pi=
            \left(%
                \begin{array}{ccccc}
                    1 & 2 & 3 & 4 & 5 \\
                    2 & 4 & 3 & 1 & 5 \\
                \end{array}%
            \right)$ will be written as $\pi=[2,4,3,1,5]$}.
        \end{rem}

        \begin{df} The \emph{symmetric group} of order $n \in \textbf{N}$
            (denoted $S_n$) is the group consisting of all the permutations
            of order $n$, with composition as the group operation.
        \end{df}

        \begin{df}\label{CoxeterA} The Coxeter generators of $S_n$ are $s_1,s_2,\ldots,s_{n-1}$
            where \\ $s_i:=[1,2,\ldots,i+1,i,\ldots,n]$.
        \end{df}
        It is a well-known fact that the \emph{symmetric group} is a \emph{Coxeter
        group} with respect to the above generating set $\{ s_i \ |\ 1 \leq i \leq n-1\}$.
        This fact gives rise to the following natural statistic of permutations
        in the symmetric group:

        \begin{df}\label{lengthA} The \emph{length} of a permutation $\pi\in S_n$
            is defined to be: \[\ell(\pi):=\min\{\ r \geq 0 \ | \  \pi=s_{i_1}
            \ldots s_{i_r} \mbox{ for some } i_1,\ldots,i_r \in [1,n] \ \}.\]
        \end{df}
        Here are other useful statistics on $S_n$ that we are going to work with:
        \begin{df} Let $\pi \in S_n$. Define the following:
            \begin{enumerate}
                \item The \emph{inversion number} of
                    $\pi$: \ \[inv(\pi):=|\{(i,j) \ | \ 1 \leq i<j\leq n, \ \pi(i)>\pi(j)\}|.\]
                    Note that $inv(\pi)=\ell(\pi)$.
                \item The \emph{descent set} of $\pi$:
                    $Des(\pi):=\{1 \leq i \leq n-1 \ | \ \pi(i)>\pi(i+1)\}$.
                \item The \emph{decent number} of $\pi$:
                    $des(\pi)=|Des(\pi)|$.
                \item The \emph{major-index} of $\pi$: $maj(\pi):=
                    \sum\limits_{i \in Des(\pi)} i$.
                \item The \emph{sign} of $\pi$:
                    $sign(\pi):=(-1)^{\ell(\pi)}$.
                \item The \emph{excedance number} of $\pi$:
                    $exc(\pi):=|\{1 \leq i \leq n \ | \ \pi(i)>i\}|$.
            \end{enumerate}
        \end{df}

        \begin{exa} Let $\pi=[2,3,1,5,4] \in S_5$. We can compute
            the above statistics on $\pi$, namely:
            \begin{eqnarray*}
                && inv(\pi)=\ell(\pi)=3, \ Des(\pi)=\{2,4\}, \ des(\pi)=2, \ maj(\pi)=6, \\
                &&  sign(\pi)=(-1)^3=-1, \mbox{ and } exc(\pi)=3.
            \end{eqnarray*}
        \end{exa}

        \begin{rem} Throughout the paper we use the following notations for a nonnegative integer $n$:
            \begin{eqnarray*}
                &&  [n]_q:=\frac{1-q^n}{1-q}, \
                    [n]_q!=[1]_q[2]_q\ldots[n]_q, \\
                &&  [n]_{\pm q}!:=[1]_q[2]_{-q}[3]_q[4]_{-q} \ldots
                    [n]_{(-1)^{n-1}q}, \mbox{ and also } \\
                &&  (a;q)_n:=\left\{%
                        \begin{array}{ll}
                            1, & \hbox{if $n=0$;} \\
                            (1-a)(1-aq)\ldots(1-aq^{n-1}), & \hbox{otherwise.} \\
                        \end{array}%
                        \right.
            \end{eqnarray*}
        \end{rem}

        MacMahon \cite{MM} was the first to find a connection between
        these statistics. He discovered that the \emph{excedance} number is equidistributed with the \emph{descent} number, and
        that the \emph{inversion} number is equidistributed with the \emph{major index}:

        \begin{thm}\label{MM}\cite{MM}
            \[\sum_{\pi \in S_n} q^{inv(\pi)}=\sum_{\pi \in S_n}
            q^{maj(\pi)}=[1]_q[2]_q[3]_q \ldots [n]_q=[n]_q!.\]
        \end{thm}

        \begin{thm}\label{desexcA}\cite{MM}
            \[\sum_{\pi \in S_n} q^{exc(\pi)}=\sum_{\pi \in S_n}
            q^{des(\pi)}.\]
        \end{thm}

        Gessel and Simion gave a similar factorial type product formula for the \emph{signed Mahonian}:
        \begin{thm}\cite[Cor.~2]{Wachs}
            \[\sum_{\pi \in S_n}
                {sign(\pi)q^{maj(\pi)}}=[n]_{\pm q}!.\] \label{GS}
        \end{thm}

        A bivariate generalization of MacMahon's Theorem~\ref{MM} was achieved during the
        1970's by Foata and Sch\"{u}tzenberger :

        \begin{thm}\label{FS1}\cite{FS}
            \[\sum_{\pi \in S_n} {q^{maj(\pi)}t^{des(\pi^{-1})}= \sum_{\pi \in
            S_n} {q^{inv(\pi)}t^{des(\pi^{-1})}}}.\]
        \end{thm}
        In the same article Foata and Sch\"{u}tzenberger
        also proved another bivariate connection between the different
        statistics:

        \begin{thm}\label{FS2}\cite{FS}
            \[\sum_{\pi \in S_n} {q^{maj(\pi^{-1})}t^{maj(\pi)}=\sum_{\pi \in
            S_n} q^{\ell(\pi)}t^{maj(\pi)}}. \]
        \end{thm}

        In 1990, during her research of the genus zeta function
        Denert, found a new statistic which was also
        Mahonian:
        \begin{df}\label{denA}\cite{Den} Let be $\pi \in S_n$, define the \emph{Denert's
            statistic} to be:
            \begin{eqnarray*}
                den(\pi)&:=&|\{1 \leq l<k \leq n \ | \ \pi(k)<\pi(l)<k\}|\\
                        &+&|\{1 \leq l<k \leq n \ | \ \pi(l)<k<\pi(k)\}|\\
                        &+&|\{1 \leq l<k \leq n \ | \ k<\pi(k)<\pi(l)\}|.
            \end{eqnarray*}
        \end{df}

        Later in the same year Foata and Zeilberger proved that
        the pair of statistics $(exc,den)$ is equidistributed with
        the pair $(des,maj)$:
        \begin{thm}\label{FZ}\cite{FZ}
            \[\sum_{\pi \in S_n} q^{exc(\pi)}t^{den(\pi)}=
                \sum_{\pi \in S_n} q^{des(\pi)}t^{maj(\pi)}.\]
        \end{thm}

    \subsection{Signed Permutations Groups}
        In this subsection we present the main definitions, notation and
        theorems for the classical Weyl groups of type B, also known as the
        \emph{hyperoctahedral groups} or the \emph{signed permutations groups}, and denoted
        $B_n$.

        \begin{df} The \emph{hyperoctahedral group} of order $n \in \textbf{N}$
            (denoted $B_n$) is the group consisting of all the bijections
            $\sigma$ of the set $[-n,n]\backslash\{0\}$ onto itself such that
            $\sigma(-a)=-\sigma(a)$ for all $a \in [-n,n]\backslash\{0\}$,
            with composition as the group operation.
        \end{df}

        \begin{rem} There are different notations for a permutation
            $\sigma \in B_n$. We will use the notation $\sigma=[\sigma(1),...,\sigma(n)]$.
        \end{rem}

        We identify $S_n$ as a subgroup of $B_n$, and $B_n$ as a subgroup
        of $S_{2n}$. As in $S_n$ we also have many different statistics; we will
        describe the main ones:

        \begin{thm} Let $\sigma \in B_n$, define the following
            statistics on $\sigma$:
            \begin{enumerate}
                \item The \emph{inversion number} of $\sigma$: \
                    $inv(\sigma):=|\{(i,j) \ | \   1 \leq i<j\leq n, \ \sigma(i)>\sigma(j) \}|.$
                \item The \emph{descent set} of $\sigma$:
                    \[Des(\sigma):=\{ 1 \leq i \leq n-1 \ | \  \sigma(i)>\sigma(i+1) \}.\]
                \item The \emph{type A descent number} of $\sigma$: \
                    $des_A(\sigma):=|Des(\sigma)|$.
                \item The \emph{type B descent number} of $\sigma$: \
                    \[des_B(\sigma):=|\{0 \leq i \leq n-1 \ | \  \sigma(i)>\sigma(i+1) \}|\mbox{, where here } \sigma(0):=0.\]
                \item The \emph{major index} of \ $\sigma$: $maj(\sigma):=
                    \sum\limits_{i \in Des(\sigma)} i.$
                \item The \emph{negative set} of $\sigma$: \ $Neg(\sigma):=\{i \in
                    [1,\ldots,n] \ | \ \sigma(i)<0 \ \}$.
                \item The \emph{negative number} of $\sigma$: \
                    $neg(\sigma):=|Neg(\sigma)|$.
                \item The \emph{negative number sum} of $\sigma$: \
                    $nsum(\sigma):=-\sum\limits_{i \in Neg(\sigma)} \sigma(i)$.
            \end{enumerate}
        \end{thm}

        \label{CoxeterB}
        It is well known (see, e.g. \cite[Proposition $8.1.3$]{BB}) that
        $B_n$ is a Coxeter group with respect to the generating set
        $\{s_0,s_1,\ldots,s_{n-1}\}$, where $s_i, \ 1 \leq i \leq n-1$, are defined as in
        $S_n$ (see~\ref{CoxeterA}), and $s_0$ is defined as: \[s_0:=[-1,2,3,\ldots,n].\] This gives
        rise to another natural statistic on $B_n$, the \emph{length}
        statistic:

        \begin{df} For all $\sigma \in B_n $ the \emph{length} of $\sigma$ is:
            \[\ell(\sigma):=\min\{r \geq 0 \ | \  \sigma=s_{i_1}s_{i_2} \ldots
            s_{i_n} \mbox{ for some } i_1,\ldots,i_r \in [0,n-1]\}.\]
        \end{df}
        There is a well-known direct combinatorial way to compute this
        statistic:

        \begin{thm}(\cite[Propositions $8.1.1$ and $8.1.2$]{BB})
            For all $\sigma \in B_n$ the length of $\sigma$ can be computed as:
            \[\ell(\sigma)=inv(\sigma)-\sum_{i \in Neg(\sigma)} \sigma(i).\]
        \end{thm}

        Using the last definition we can define another natural
        statistic on $B_n$, the \emph{sign} statistic:
        \begin{df} For all $\sigma \in B_n$ the \emph{sign} of $\sigma$ is:
            \[sign(\sigma):=(-1)^{\ell(\sigma)}.\]
        \end{df}
        The generating function of \emph{length} is also called the
        Poincar\'{e} polynomial and can presented
        in the following manner:

        \begin{thm}\cite[\S 3.15]{Hu}
           \[\sum_{\sigma \in B_n} q^{\ell(\sigma)}=[2]_q[4]_q \ldots [2n]_q=\prod_{i=1}^{n} [2i]_q.\]
        \end{thm}
        Recently, Adin and Roichman generalized MacMahon's result Theorem~\ref{MM} to $B_n$, by
        introducing a new Mahonian statistic, the \emph{flag major index}:
        \begin{df}\label{fmaj}\cite{AR} The \emph{flag major index} of $\sigma \in B_n$ is defined
            as:\[\mbox{flag-major}(\sigma):=2maj(\sigma)+neg(\sigma),\] where
            $maj(\sigma)$ is calculated with respect to the linear order
            \[-1<-2<\ldots<-n<1<2<\ldots<n.\]
        \end{df}

        \begin{thm}\label{fmajlen}\cite[\S 2]{AR}
            \[\sum_{\sigma \in B_n} q^{\ell(\sigma)}=
            \sum_{\sigma \in B_n} q^{flag-major(\sigma)}=[2]_q[4]_q \ldots [2n]_q.\]
        \end{thm}

        \begin{rem}The previous result still holds if $maj(\sigma)$ is calculated
            with respect to the \emph{natural order}
            $-n<-(n-1)<...<-2<-1<1<2<...<n-1<n,$ see also \cite{AR}.
        \end{rem}

        Adin, Brenti and Roichman introduced another statistic
        which was also Mahonian, the \emph{nmaj} statistic:
        \begin{df}\label{nmaj}\cite[\S 3.2]{ABR} Let $\sigma \in B_n$ then the \emph{negative
            major index} is defined as:
            \[nmaj(\sigma):=maj(\sigma)-\sum_{i \in Neg(\sigma)} \sigma(i)=maj(\sigma)+nsum(\sigma).\]
        \end{df}

        \begin{thm}\label{lenmaj}\cite{ABR}
            \[\sum_{\sigma \in B_n} q^{\ell(\sigma)}=\sum_{\sigma \in B_n} q^{nmaj(\sigma)}.\]
        \end{thm}

        In the same article \cite{ABR} they also defined a new descent
        multiset and new descent statistics, and found a new Euler-Mahonian bivariate
        distribution for these statistics:
        \begin{df}\cite[\S 3.1 and \S 4.2]{ABR} Let $\sigma \in B_n$ define:
            \begin{enumerate}
                \item  The \emph{negative descent multiset} of
                    $\sigma$: \[NDes(\sigma):=Des(\sigma)\bigcup \{-\sigma(i)\ |\
                    i \in Neg(\sigma)\},\] where $\bigcup$ stands
                    for multiset union.
                \item The \emph{negative descent} statistic of
                    $\sigma$: $ndes(\sigma):=|NDes(\sigma)|$.
                \item The \emph{flag-descent} number of $\sigma$:
                    $fdes(\sigma):=des_A(\sigma)+des_B(\sigma)=2des_A(\sigma)+\varepsilon(\sigma)$,
                    where
                    \begin{eqnarray*}
                        \varepsilon(\sigma):=
                        \left\{%
                            \begin{array}{ll}
                                1, & \hbox{if $\sigma(1)<0$;} \\
                                0, & \hbox{otherwise.} \\
                            \end{array}%
                        \right.
                    \end{eqnarray*}
            \end{enumerate}
        \end{df}

        \begin{thm}\label{fndesmaj}\cite[\S 4.3]{ABR}
            \[\sum_{\sigma \in B_n} t^{ndes(\sigma)}q^{nmaj(\sigma)}=
                \sum_{\sigma \in B_n} t^{fdes(\sigma)}q^{flag-major(\sigma)}.\]
        \end{thm}
        In their article from 2005 Adin, Gessel,
        and Roichman gave a type B analogue to the Gessel-Simion
        Theorem(e.g. \cite[Cor.~2]{Wachs}):
        \begin{thm}\cite[\S 5.1]{ARG}
            \[ \sum_{\sigma \in B_n}
            sign(\sigma)q^{flag-major(\sigma)}=[2]_{-q}[4]_q \ldots
            [2n]_{(-1)^nq}.\]
            Where \emph{flag major index} computed with respect to the
            linear order: \[-1<-2<\ldots<-n<1<2<\ldots<n.\]
        \end{thm}
        The group $B_n$ has some well known subgroups
        \cite[\S7]{ARG}, in this paper we will only work with the subgroup of
        \emph{elements with even neg}:

        \begin{df}\label{Dn}\label{wsubgroups} The subgroup of \emph{elements with even neg} in $B_n$ (denoted
        $D_n$) is defined as:
            \[D_{n}:=\{\sigma \in B_n \ | \ neg(\sigma)\equiv 0 \ mod \ 2 \}.\]
        \end{df}
    \medskip

    \subsection{Colored Permutation Groups}
        In this subsection we will present the definitions,
        notations, and basic theorems for \emph{the colored permutation groups} that we are
        going to work with throughout this paper. We start with the definition of the \emph{colored permutation
        groups}.

        \begin{df} Let $r,n$ be positive integers. The group of
            colored permutations of $n$ digits and $r$ colors, denoted
            $G_{r,n}$, is the wreath product $\ G_{r,n}=C_r \wr S_n$,
            consisting of all the pairs $(z,\pi)$ where $z$ is a n-tuple of integers between
            $0$ and $r-1$, and $\pi \in S_n$. The group operation defined as follows:
            for $z=(z_1,\ldots,z_n), \ z'=(z_1',\ldots,z_n'), \ \pi, \ \pi' \in S_n$,
            \[(z,\pi)\cdot(z',\pi')=((z_1+z'_{\pi^{-1}(1)},\ldots,z_n+z'_{\pi^{-1}(n)}),\pi
            \circ\pi')\]
            (where $+$ is taken $mod \ r$).
        \end{df}

        \begin{rem} Let $r$ and $n$ be nonnegative integers. We will use the notation \\ $g,\ \bar{g} \in G_{r,n}$,
            where $g=(z,\pi)$ and $\bar{g}=(z,\pi^{-1})$.
        \end{rem}

        \begin{nt} Notice that for  $r=1$, $G_{1,n}$ isomorphic to the \emph{symmetric group}
            $S_n$, and for $r=2$, $G_{2,n}$ isomorphic to the \emph{signed
            permutation} group $B_n$.
        \end{nt}
        In the following definitions we assume that the alphabet:
        \[\{1,\ldots,n,1^{[1]},\ldots,n^{[1]},\ldots,1^{[r-1]},\ldots,n^{[r-1]}\}\]
        has some predefined linear order. As in the cases of the \emph{symmetric
        group}, and the \emph{signed permutation group} we can also
        define statistics on the \emph{colored permutation groups}:

        \begin{df} Let be $g=(z,\pi) \in G_{r,n}$. Define:
            \begin{enumerate}
                \item The \emph{inversion number} of $g$:
                    \[inv(g):=|\{(i,j) \ | \ 1 \leq i <j \leq n,
                    \ \pi(i)^{[z_i]}>\pi(j)^{[z_j]} \}|.\]
                \item The \emph{Descent set} of $g$:
                    $Des(g):=\{\ i \ |\ 1 \leq i \leq n-1, \
                    \pi(i)^{[z_i]}>\pi(i+1)^{[z_{i+1}]} \ \}$.
                \item The \emph{descent number} of $g$:
                    $des(g):=|Des(g)|.$
                \item The \emph{major index} of $g$:
                    $maj(g):=\sum\limits_{i \in Des(g)}i$.
                \item The \emph{Negative set} of $g$:
                    $Neg(g):=\{\ i \ | \ z_i \neq 0 \ \}$.
                \item The \emph{negative number} of $g$:
                    $neg(g):=|Neg(g)|$.
                \item The \emph{color sum} of $g$:
                    $csum(g):=\sum\limits_{z_i \neq 0} z_i.$
                \item The \emph{negative color sum} of $g$:
                    $ncsum(g):=\sum\limits_{z_i \neq 0}
                    (\pi(i)-1).$
            \end{enumerate}
        \end{df}

        \label{rlength}

        The \emph{colored permutation group} $G_{r,n}$ has a natural
        \emph{length} statistic, defined with respect to generating set
        $\{s_0,s_1,\ldots,s_{n-1}\}$, where $s_i, \ i\geq 1$ is
        defined as\\ $s_i:=((0,\ldots,0),[1,\ldots,i+1,i,\ldots,n])$ and $s_0:=((1,0,\ldots,0),id)$,
        see also \cite{ST} and \cite{BE}. There is a direct algebraic formula to
        calculate the \emph{length} of $g \in G_{r,n}$, namely:
        \begin{thm}\label{rlen}\cite{BE}
            \[\ell(g):=inv(g)+\sum\limits_{z_i \neq 0} (\pi(i)-1)+\sum\limits_{i=1}^n
            z_i,\] where $inv(g)$ is calculate
            according to the alphabet linear order:
                \[n^{[r-1]}<\ldots<n^{[1]}<\ldots<1^{[r-1]}<\ldots<1^{[1]}<1<\ldots<n.\]
        \end{thm}
        Adin and Roichman found a \emph{flag major index}
        statistic on \emph{colored permutation group}:
        \begin{df}\label{rfmaj} \cite[\S 3.1]{AR} Let be $g \in G_{r,n}$,
            define:
            \[flag-major(g):= r \cdot
            maj(g)+\sum_{z_i \neq 0} z_i.\]
        \end{df}

        In the case of \emph{colored permutation groups} we also have new Eulerian and Mahonian statistics which have
        been found in the passing year by Bagno:

        \begin{df} Let be $g \in G_{r,n}$, define:
            \begin{enumerate}
                \item $ldes(g):=des(g)+\sum\limits_{i=1}^n z_i.$
                \item $lmaj(g):=maj(g)+\sum\limits_{z_i
                    \neq 0} (\pi(i)-1)+\sum\limits_{i=1}^n z_i.$
            \end{enumerate}
        \end{df}

        \begin{thm}\cite[\S 5]{BE}
            \[\sum_{g \in G_{r,n}} q^{lmaj(g)}=
                \sum_{g \in G_{r,n}}
                q^{\ell(g)}=[n]_q!\prod_{i=1}^n(1+q^i[r-1]_q).\]
        \end{thm}

        Bagno also found an interesting statistic over $G_{r,n}$, the \emph{nmaj} statistic
        (see definition \cite[\S 6]{BE}), which generalize the type B \emph{nmaj} statistic
        (see~\ref{nmaj}), and proved the following equality:
        \begin{thm}\label{rnmaj}\cite[\S 7]{BE}
            \[\sum_{g \in G_{r,n}}
            q^{nmaj(g)}=\sum_{g \in G_{r,n}}
            q^{flag-major(g)}=\prod_{i=1}^n [ri]_q.\]
        \end{thm}

\section{Main Results}
    In this section we present the main results obtained in the paper according to their order of appearance:

    \begin{rem} Throughout this section we use the notation $stat_{K_n}(g)$
        when \it{stat} is a statistic on $G_{r,n}$ and $K_n$ is a predefined
        linear order on the alphabet of $G_{r,n}$. The notation $stat_{K_n}(g)$
        means that we calculate the statistic \emph{stat} according to
        the linear order $K_n$.
    \end{rem}
    We begin Section 4 with Proposition~\ref{fmajorder}, which states that the \emph{flag major index}
    is equidistributed with the \emph{length} over the \emph{signed
    permutation group} for any linear order:
    \begin{pro} (Proposition~\ref{fmajorder}) For any linear order $K_n$:
        \[\sum_{\sigma \in B_n} q^{flag-major(\sigma)}=\sum_{\sigma \in B_n}
        q^{flag-major_{K_n}(\sigma)}.\]
    \end{pro}
    The main results of Section 5 calculate generating functions according to \emph{natural order}:
    \[\ N: -n<-(n-1)<\ldots<-1<1<\ldots<n-1<n,\] over the \emph{signed permutations groups}.
    We begin with Theorem~\ref{BSF}, which gives the generating
    function of the \emph{signed Mahonian}:
    \begin{thm} (Theorem~\ref{BSF})
        \[\sum\limits_{\sigma \in B_{n}} sign(\sigma) q^{flag-major_N(\sigma)}=
            (q;-1)_n[n]_{\pm q^2}!.\]
    \end{thm}
    The next result of this type:
    \begin{thm}(Theorem~\ref{signmaj})
        \[ \sum_{\sigma \in B_n} sign(\sigma)q^{nmaj_N(\sigma)}=(q;-q)_n[n]_{\pm q}!.\]
    \end{thm}
    In Theorem~\ref{fnmaj}, we calculate an interesting bivariate
    generating function of two Mahonian statistics:
    \begin{thm}(Theorem~\ref{fnmaj})
        \[\sum_{\sigma \in B_n} q^{flag-major_N(\sigma)}t^{nmaj_N(\sigma)}=\prod_{i=1}^n(1+qt^i)[n]_{q^2t}!.\]
    \end{thm}
    In Theorem~\ref{fmajlenN} we calculate another interesting bivariate generating function of \emph{length} and
    \emph{flag major index}:
    \begin{thm}(Theorem~\ref{fmajlenN})
        \[\sum_{\sigma \in B_n} q^{flag-major_N(\sigma)}t^{\ell(\sigma)}=
            \prod_{i=1}^n (1+qt^i)A_n(q^2,t),\] where
            $A_n(q,t)=\sum\limits_{\pi \in S_n} q^{maj(\pi)}t^{\ell(\pi)}=
            \sum\limits_{\pi \in S_n} q^{maj(\pi)}t^{inv(\pi)}.$
    \end{thm}
    We finish reviewing Section 5 with Theorem~\ref{DSF}, which presents a calculation of the generating function of
    \emph{sign} and \emph{flag major index} over $D_n$:
    \begin{thm}(Theorem~\ref{DSF})
        \begin{eqnarray*}
            \sum\limits_{\sigma \in D_n} sign(\sigma)q^{flag-major_N(\sigma)}=
            (1-q^2)^{\lfloor\frac{n}{2}\rfloor}[n]_{\pm q^2}!.
        \end{eqnarray*}
    \end{thm}
    According to the previous theorems we can conclude very interesting result, which is presented in Corollary~\ref{DSFC}:
    \begin{cor}(Corollary~\ref{DSFC})
            \[ \sum\limits_{\sigma \in D_{2n}} sign(\sigma)q^{flag-major_N(\sigma)}=
            \sum\limits_{\sigma \in B_{2n}} sign(\sigma)q^{flag-major_N(\sigma)}.\]
    \end{cor}
    We move on to Section 6, in which we calculate generating
    functions over the \emph{colored permutation groups} according
    to the \emph{friends order}:
    \[F: \ 1^{[r-1]}<\ldots<1<2^{[r-1]}<\ldots<2<\ldots<n^{[r-1]}<\ldots<n.\]
    The first main result in this section is Lemma~\ref{tpc},
    which is used in almost all the proofs in this
    section:
    \begin{lem}(Lemma~\ref{tpc}) Let $i \in [1,k]$ and let $ostat_i: \ S_n \to \mathbf{Z}, \ cstat_i: \
        \mathbf{Z}_r^n \to \mathbf{Z}$ statistics on
        $S_n$ and $\mathbf{Z}_r^n$ respectively; define $stat_i: G_{r,n} \to
        \mathbf{Z}$ by $stat_i(z,\pi):=ostat_i(\pi)+cstat(z)$.
        Then:
        \[\sum_{(z,\pi) \in G_{r,n}} q_1^{stat_1(z,\pi)}
            \ldots q_k^{stat_k(z,\pi)}=
            \sum_{(z,id) \in Color^n_r} q_1^{cstat_1(z)}
            \ldots q_k^{cstat_k(z)} \sum_{\pi \in S_n} q_1^{ostat_1(\pi)}
            \ldots q_k^{ostat_k(\pi)},\]
            where the set $Color^n_r$ is defined in~\ref{COLORN}.
    \end{lem}
    We continue on, define a new \emph{friends color-sign} statistic (see Definition~\ref{fcsign}), and
    calculate the bivariate generating function of \emph{friends color-sign} and \emph{flag major index} over
    the \emph{colored permutation groups}. The result of this calculation
    is presented in Theorem~\ref{csignfmaj}:
    \begin{thm}(Theorem~\ref{csignfmaj})
           \[\sum_{g \in G_{r,n}} csign_F(g)q^{flag-major_F(g)}=
           [r]_{-q}^n[n]_{\pm rq}!.\]
    \end{thm}
    We continue this section with the definitions of two new
    statistics over: the \emph{colored permutations groups}: the \emph{r-color excedance number}, denoted $exc_r$ (see definition~\ref{excr}),
    and the \emph{r-color Denert's} statistic, denoted $den_r$ (see
    Definition~\ref{denr}). By using these new statistics we give in Theorem~\ref{FZR}
    a generalization of the Foata-Zeilberger Theorem~\ref{FZ}:
    \begin{thm} (Theorem~\ref{FZR}) Let $r,n$ be positive integers. Then:
        \[\sum_{g \in G_{r,n}} q^{den_r(g)}t^{exc_r(g)}=\sum_{g \in
        G_{r,n}} q^{flag-major_F(g)}t^{ldes_F(g)}.\]
    \end{thm}
    We finish the paper in Section 7, where we define new statistic, the \emph{flag-inversion} statistic (see Definition~\ref{finv})
    which is equidistributed with \emph{flag major index}, and by
    using the \emph{flag-inversion} statistic we give in Theorem~\ref{FSR}, and in Theorem~\ref{GFS2} another
    generalizations to the Foata-Sch\"utzenberger Theorem~\ref{FS1}, and Theorem~\ref{FS2}:
     \begin{thm}(Theorem~\ref{FSR})
        \[ \sum_{g \in G_{r,n}} q^{flag-major_F(g)}
            t^{ldes_F(\bar{g})}=\sum_{g \in G_{r,n}} q^{finv_F(g)}
            t^{ldes_F(\bar{g})}.\]
    \end{thm}
    \begin{thm}(Theorem~\ref{GFS2})
    \[\sum_{\sigma \in B_n} q^{flag-major_F(\sigma^{-1})}t^{flag-major_F(\sigma)}=
        \sum_{\sigma \in B_n} q^{finv_F(\sigma)}t^{flag-major_F(\sigma)}.\]
    \end{thm}

\section{Flag Major Index  With Different Orders}
    In this section study statistics on the \emph{signed permutation group}
    with respect to different linear orders. First, we give necessary definitions and
    introduce certain bijections (see Definition~\ref{ibij}). Using
    these definitions we prove that the \emph{flag major
    index} is equidistributed with the \emph{length} statistic for any linear
    order; however, this result does not hold in the case of \emph{signed
    Mahonian} distribution (it is easy to check this fact in the case of $n=2$),
    we finish this section indicating a large collection of linear orders
    which are \emph{sign} and \emph{flag major index} equidistributed with
    each other.

    \begin{df}\label{loeders} A \emph{linear order of length} \it{n}, denoted $K_n$, is a
        bijection \[K_n:[1,2n] \to[-n,n]\backslash
        \{0\}.\]
        For convenience we write a linear order in the following
        form:
        \[K_n: \ K_n(1)<K_n(2)<\ldots<K_n(2n-1)<K_n(2n).\]
    \end{df}

    We can calculate permutation statistics according to a linear
    order $K_n$, we use the following notation: $maj_{K_n}(\sigma),\ des_{K_n}(\sigma),\
    flag-major_{K_n}(\sigma),\ nmaj_{K_n}(\sigma)$ etc, to indicate that the corresponding statistic
    is calculated with respect to the linear order $K_n$. We also use the
    notation: $m>_{K_n}l,$ to indicate, that according to the linear order $K_n$ 'm' is larger than 'l',
    i.e. that $m=K_n(s)$, \ $l=K_n(r)$, and $s>r$.

    \begin{exa} Let $K_n$ be a linear order and let $\sigma \in B_n$. Then:
        \[maj_{K_n}(\sigma):=\sum_{\sigma(i)>_{K_n}\sigma(i+1)} {i}. \]
    \end{exa}

    \begin{nt} Notice that for any linear order $K_n$, and for any
        $\sigma \in B_n, \ neg(\sigma)=neg_{K_n}(\sigma)$. This also
        applies to the \emph{length} statistic, because it is defined
        with respect to the Coxeter generators, which do not depend on the choice
        of linear order.
    \end{nt}

    \begin{df} Let $K_n: \ K_n(1)<K_n(2)<\ldots<K_n(2n-1)<K_n(2n),$ be a linear order and let $1 \leq j\leq 2n-1$.
    Define $K_{n,j}$ to be the following linear order:
        \[K_{n,j}: \ K_n(1)<K_n(2)< \ldots<K_n(j+1)<K_n(j)< \ldots <K_n(2n-1)<K_n(2n).\]
    \end{df}

    \begin{df}\label{ibij} Let $K_n$ be a linear order, $j \in [1,2n-1]$, and $s \in [1,n]$. Define
        the function $\psi^{j}(\sigma):B_n \rightarrow B_n$ as:
        \begin{center}
            $\psi^{j}(\sigma)(s):=
            \left\{%
                   \begin{array}{ll}
                        K_n(j), & \hbox{$\sigma(s)=K_n(j+1),$}\\
                        \    &  \hbox{$\exists p \in [1,n], \ \sigma(p)=K_n(j)$ }, \\
                        K_n(j+1), & \hbox{$\sigma(s)=K_n(j),$}\\
                        \   &   \hbox{$\exists p \in [1,n], \ \sigma(p)=K_n(j+1),$}\\
                        \sigma(s), & \hbox{Otherwise.}\\
                    \end{array}%
            \right.$
        \end{center}
    \end{df}

    \begin{nt} By the definition of $\psi^j$ we can conclude that $\psi^{j}(\psi^j(\sigma))=\sigma$.
    \end{nt}

    \begin{lem} Let $K_n$ be a linear order and $j \in [1,2n-1]$, then $\psi^{j}:B_n \rightarrow B_n$
        is a bijection, and \[Des_{K_{n}}(\sigma)=Des_{K_{n,j}}(\psi^j(\sigma)), \ \forall \sigma \in
        B_n.\]
    \end{lem}

    \noindent{\bf Proof.} Let $\sigma_1,\ \sigma_2 \in B_n$, we need to
    prove that: $\sigma=\psi^j(\sigma_1)=\psi^j(\sigma_2)\Rightarrow
    \sigma_1=\sigma_2$. We divide our proof into two cases:
    \begin{enumerate}
        \item  If $\exists p,m \in [1,n]$, such that $\sigma(p)=K_n(j), \
            \sigma(m)=K_n(j+1)$, then according to the definition of
            $\psi^j$ we conclude that
            $\sigma=\psi^j(\sigma_1)=\psi^j(\sigma_2) \Rightarrow \sigma_1(s)=\sigma_2(s), \forall s \neq
            p,m$, and $\sigma(p)=\sigma_1(m)=\sigma_2(m), \
            \sigma(m)=\sigma_1(p)=\sigma_2(p)$, therefore we
            conclude that: $\sigma_1=\sigma_2$.
        \item Otherwise, according to the definition of $\psi^j$ we
             get $\psi^j(\sigma_1)=\sigma_1=\psi^j(\sigma_2)=\sigma_2 \Rightarrow
             \sigma_1=\sigma_2$.
    \end{enumerate}
    We have proved that $\psi^j$ is a bijection. Now we
    need to prove that: $Des_{K_n}(\sigma)=Des_{K_{n,j}}(\psi^j(\sigma))$.
    We will do so by checking all the possible cases of
    $\psi^j(\sigma)(i)$ and
    $\psi^j(\sigma)(i+1)$, and proving that for each case $\sigma(i)>_{K_n} \sigma(i+1)
    \Leftrightarrow \psi^j(\sigma)(i)>_{K_{n,j}} \psi^j(\sigma)(i+1)$:
    \begin{enumerate}
        \item \(\sigma(i)=K_n(s), \ s \in [1,2n], \ s \neq j,j+1\), therefore by the definition
            of $\psi^j$:\\ $\psi^j(\sigma)(i)=\sigma(i)=K_n(s)=K_{n,j}(s)$, and again
            we divide into cases:
            \begin{enumerate}
                \item $\sigma(i+1)=K_n(t), \ t \in [1,2n], \ t \neq j,j+1$
                    therefore; $\psi^j(\sigma)(i+1)=\sigma(i+1)=K_n(t)=K_{n,j}(t),$ and now we can
                    conclude that:
                        \begin{eqnarray*}
                            &&  \sigma(i)=K_n(s)>_{K_n} K_n(t)=\sigma(i+1) \Leftrightarrow s>t
                                \Leftrightarrow \sigma(i)>_{K_{n,j}} \sigma(i+1) \Leftrightarrow \\
                            &&  \psi^j(\sigma)(i)>_{K_{n,j}} \psi^j(\sigma)(i+1).
                        \end{eqnarray*}
                \item $\sigma(i+1)=K_n(j)$, again we divide into
                    two cases according to the behavior of $\psi^j$:
                    \begin{enumerate}
                        \item $\exists p \in [1,n], \ \sigma(p)=K_{n}(j+1)$, in
                            this case:
                            \begin{eqnarray*}
                                &&  \psi^j(\sigma)(i+1)=K_n(j+1)=K_{n,j}(j), \mbox{ and } \\
                                &&  \psi^j(\sigma)(i)=K_n(s)=K_{n,j}(s),\mbox{ therefore;}\\
                                &&  \sigma(i)>_{K_n} \sigma(i+1) \Leftrightarrow s>j \Leftrightarrow
                                    K_{n,j}(s)>K_{n,j}(j) \Leftrightarrow\\
                                &&  \psi^j(\sigma)(i)>_{K_{n,j}} \psi^j(\sigma)(i+1).
                            \end{eqnarray*}

                        \item $\forall p \in [1,n], \
                            \sigma(p)\neq K_{n}(j+1),$ in this case $\psi^j(\sigma)=\sigma$
                            \begin{eqnarray*}
                                &&  \psi^j(\sigma)(i+1)=K_{n}(j)=K_{n,j}(j+1), \mbox{ therefore;} \\
                                &&  \sigma(i)>_{K_n} \sigma(i+1) \Leftrightarrow s>j \Leftrightarrow_{s \neq j,j+1}\\
                                &&  s>j+1 \Leftrightarrow
                                    \psi^j(\sigma)(i)=K_{n,j}(s)>_{K_{n,j}}
                                    K_{n,j}(j+1)=\psi^j(\sigma)(i+1).
                            \end{eqnarray*}
                \end{enumerate}

            \item  $\sigma(i+1)=K_n(j+1),$ this case is almost identical to the previous case (b) with minor
                    changes.
            \end{enumerate}

        \item $\sigma(i)=K_n(j)$, as in the previous case, we will check all the
            possible options:
            \begin{enumerate}
                \item $\sigma(i+1)=K_n(t), \  t \in [1,2n], \ t \neq
                    j+1$, and we divide into cases according to the behavior of $\psi^j$:
                    \begin{enumerate}
                        \item $\exists p\in [1,n] \ \sigma(p)=K_n(j+1),$
                        then:
                            \begin{eqnarray*}
                                &&  \psi^j(\sigma)(i)=K_n(j+1)=K_{n,j}(j),\\
                                &&  \psi^j(\sigma)(i+1)=K_n(t)=K_{n,j}(t), \mbox{ therefore;}\\
                                &&  \sigma(i)>_{K_n} \sigma(i+1) \Leftrightarrow
                                    j>t \Leftrightarrow \\
                                &&  \psi^j(\sigma)(i)=K_{n,j}(j)>_{K_{n,j}}
                                    K_{n,j}(t)=\psi^j(\sigma)(i+1).
                            \end{eqnarray*}

                        \item $\forall p \in [1,n] \ \sigma(p)\neq K_n(j+1)$, then:
                            \begin{eqnarray*}
                                &&  \psi^j(\sigma)(i)=K_n(j)=K_{n,j}(j+1),\\
                                &&  \psi^j(\sigma)(i+1)=K_n(t)=K_{n,j}(t), \mbox{ therefore;}\\
                                &&  \sigma(i)>_{K_n} \sigma(i+1) \Leftrightarrow  j>t {\Leftrightarrow} \\
                                &&   j+1>t \Leftrightarrow
                                    \psi^j(\sigma)(i)=K_{n,j}(j+1)>_{K_{n,j}}
                                    K_{n,j}(t)=\\
                                &&  \psi^j(\sigma)(i+1).
                            \end{eqnarray*}
                    \end{enumerate}
                \item $\sigma(i+1)=K_n(j+1)=K_{n,j}(j)$, then:
                    \begin{eqnarray*}
                        &&  \psi^j(\sigma){i}=K_n(j+1)=K_{n,j}(j), \\
                        &&  \psi^j(\sigma)(i+1)=K_n(j)=K_{n,j}(j+1), \mbox{ therefore; }\\
                        &&  \sigma(i+1)=K_n(j+1)>_{K_n} \sigma(i)=K_n(j) \Leftrightarrow \\
                        &&  \psi^j(\sigma)(i+1)=K_{n,j}(j+1)>_{K_{n,j}} \psi^j(\sigma)(i)=K_{n,j}(j).
                    \end{eqnarray*}
            \end{enumerate}
        \item $\sigma(i)=K_n(j+1)$. This case is similar to case 2 above.
    \end{enumerate}
    We got that:
    \[\forall \sigma \in B_n, \ \forall i \in [1,n-1], \ \sigma(i)>_{K_n} \sigma(i+1) \Leftrightarrow
    \psi^j(\sigma)(i)>_{K_{n,j}}\psi^j(\sigma)(i+1).\]
    Now we can conclude that $\forall \sigma, \ \in B_n \ Des_{K_n}(\sigma)=Des_{K_{n,j}}(\psi^j(\sigma))$  $\quad \bullet$

    \begin{cor} Let $K_n$ be a linear order and $1 \leq j\leq 2n-1$ then:
        \[maj_{K_{n,j}}(\psi^j(\sigma))=maj_{K_n}(\sigma).\]
    \end{cor}

    \begin{cor} Let $K_n$ be linear order  and $1 \leq j\leq 2n-1$ then:
        \[ \sum_{\sigma \in B_n} q^{maj_{K_n}(\sigma)}=\sum_{\sigma \in B_n}
        q^{maj_{K_{n,j}}(\sigma)}.\]
    \end{cor}

    \begin{cor}\label{eqmaj} Let $K_n$ be linear order then:
        \[\sum_{\sigma \in B_n} q^{maj(\sigma)}=\sum_{\sigma \in B_n}
        q^{maj_{K_n}(\sigma)}.\]
    \end{cor}

    \noindent{\bf Proof.} We can
    transfer the given linear order $K_n$ to the natural order: $N: \
    -n<-(n-1)<\ldots<-1<1<\ldots<n-1<n$, by a finite number of adjacent transpositions.

    \begin{rem} According to the above proof for each linear order $K_n$ there exists a bijection
        $\phi:B_n \to B_n$ that, $ maj(\phi(\sigma))=maj_{K_n}(\sigma), \ \forall  \sigma$, and
        $\phi=\psi^{r_m} \circ \psi^{r_{m-1}}\circ \ldots \circ \psi^{r_1}$ where
        $\{r_i\}_{i=1}^{m}$ are positive integers.
    \end{rem}
    Now we can prove that for any linear order $K_n$ the
    \emph{flag major index} equidistributed with $length$
    statistic over the \emph{the signed permutations groups}:

    \begin{pro}\label{fmajorder} Let $K_n$ be a linear order then:
        \[\sum_{\sigma \in B_n} q^{flag-major(\sigma)}=\sum_{\sigma \in B_n}
         q^{flag-major_{K_n}(\sigma)}.\]
    \end{pro}

    \noindent{\bf Proof.} Use the bijection $\phi: B_n \to B_n$,
    and note that \[neg(\phi(\sigma)=neg(\sigma), \ \forall \sigma \in B_n.\]
    \begin{eqnarray*}
        \sum\limits_{\sigma \in B_n} q^{flag-major_{K_n}(\sigma)}&=&
            \sum\limits_{\sigma \in B_n} {q^{2maj_{K_n}(\sigma)+neg_{K_n}(\sigma)}}\\
        &=& \sum\limits_{\sigma \in B_n} {q^{2maj_{K_n}(\phi(\sigma))+neg(\phi(\sigma))}}\\
        &=& \sum\limits_{\sigma \in B_n} {q^{2maj(\sigma)+neg(\sigma)}}\\
        &=& \sum\limits_{\sigma \in B_n} {q^{flag-major(\sigma)}} \quad \bullet
    \end{eqnarray*}

    \begin{lem} Let $K_n$ be a linear order, such that, for some  $j\in
        [1,2n-1], \ K_n(j)+K_n(j+1)=0$. Then:
        $\psi^j(\sigma)=\sigma$.
    \end{lem}

    \noindent{\bf Proof.} Let be $\sigma \in B_n$, by the definition of
    $B_n$, we know that $\exists i \in [1,n]  \ \sigma(i)=\pm K_n(j)$,
    therefore; $\forall s \in [1,n] \ \sigma(s) \neq -\sigma(i)$,
    and now according to the definition of $\psi^j$ we conclude that:
    $\psi^j(\sigma)=\sigma \quad \bullet$

    \begin{cor} Let $K_n$ be a linear order, such that, for some  $j\in
        [1,2n-1], \ K_n(j)+K_n(j+1)=0$. Then:
        \[\sum_{\sigma \in B_n}{sign(\sigma)q^{flag-major_{K_n}(\sigma)}}=
        \sum_{\sigma \in B_n}{sign(\sigma)q^{flag-major_{K_{n,j}}(\sigma)}}.\]
    \end{cor}

    \noindent{\bf Proof.}
    \begin{eqnarray*}
            \sum_{\sigma \in B_n} sign(\sigma) q^{flag-major_{K_n}(\sigma)}&=&
            \sum_{\sigma \in B_n} sign(\psi^j(\sigma)) q^{flag-major_{K_n}(\psi^j(\sigma))}\\
        &=& \sum_{\sigma \in B_n} sign(\sigma) q^{flag-major_{K_{n,j}}(\sigma)}\\
        &=& \sum_{\sigma \in B_n} sign(\sigma) q^{flag-major_{K_{n,j}}(\sigma)} \quad\bullet
    \end{eqnarray*}

    \begin{rem} Let $K_n$ be linear order, and let $\{i_s\}_{s=1}^{s=m}$
        be integers in $[1,2n-1]$, such that $K_n(i_s)+K_n(i_{s+1})=0 \ \forall s$,
        then it easy to prove by induction that for any linear order
        $K'_n=K_{n,i_{j_1},i_{j_2}, \ldots ,i_{j_w}}, \ w \in \textbf{N}, \ i_{j_k} \in \{ i_s \}_{s=1}^{s=m},
        \ k \in [1,w] $, exist the following equality:
        \[\sum_{\sigma \in B_n}
        {sign(\sigma)q^{flag-major_{K_n}(\sigma)}}=\sum_{\sigma \in B_n}
        {sign(\sigma)q^{flag-major_{K'_n}(\sigma)}}.\]
    \end{rem}

    \begin{exa} For example consider the linear order:
        \[F:\ -1<1<-2<2<\ldots<-n<n,\] then we can create new linear order
        $F^{2i-1}=Fs_{2i-1}$ by swapping $F(2i-1)$ with $F(2i)$ for all
        $i \in [1,n]$. And according to the last corollary we get
        that: \[\sum_{\sigma \in B_n}
        {sign(\sigma)q^{flag-major_{F}(\sigma)}}=\sum_{\sigma \in B_n}
        {sign(\sigma)q^{flag-major_{F^{2i-1}}(\sigma)}}.\] According to the last
        remark we can generate from the linear order $F$, $2^n$ linear orders which also
        satisfy the last equality.
    \end{exa}

\section{Permutation Statistics in the Natural Order}
    In this section we will calculate generating functions of
    different statistics according to the \emph{natural order}:
    \[N: \ -n<-(n-1)<\ldots<-1<1<\ldots<n-1<n.\]
    Most of the results of this section have been achieved by reducing the summation from
    $B_n$ to $U_n$ (see also \cite[\S 7]{ARG}), and calculating the generating function
    over $U_n$ by recursion to $U_{n-1}$. Before we begin with our calculation and theorems, we need some helpful
    tools that will be presented in the following section.

    \medskip

    \subsection{Preliminaries} \label{dun}
        \begin{df} Define the following set:
            \[U_n:=\{\tau \in B_n \ | \ \tau(1)<\tau(2)<\ldots<\tau(n-1)<\tau(n)\}.\]
        \end{df}
        There are several facts (see also \cite{ABR},\cite{ARG}) about the set $U_n$ that can be directly concluded
        from the definition of $U_n$, namely:
        each $\sigma \in B_n$ has a unique representation as:
        \[\sigma=\tau\pi \ (\tau \in U_n \mbox{ , and } \pi \in S_n).\]
        The following properties are clear:
        \begin{fac}
            \begin{enumerate}
                \item $maj_N(\sigma)=maj_N(\pi).$
                \item $inv_N(\sigma)=inv_N(\pi).$
                \item $neg(\sigma)=neg(\tau).$
                \item $nsum(\sigma)=nsum(\tau).$
                \item $sign(\sigma)=sign(\tau)sign(\pi)$.
            \end{enumerate}
        \end{fac}
        We continue with our definitions.

        \begin{df} Define the following subsets of $U_n$:
            \begin{enumerate}
                \item $U_{n1}:=\{\tau \in U_n \ | \  \tau(1)=-n  \}$.
                \item $U_{n2}:=\{\tau \in U_n\ |\  \tau(n)=n \}$.
            \end{enumerate}
        \end{df}

        \begin{nt}
            $U_n=U_{n1}\uplus U_{n2}$, where $\uplus$ stands for disjoint union.
        \end{nt}
        We also define two bijections from $U_{n-1}$ one onto $U_{n1}$,
        and one onto $U_{n2}$:
        \begin{df} For $i \in {1,2}$, define $\varphi_{ni}: U_{n-1} \to U_{ni}$ by:
            \begin{enumerate}
                \item $\varphi_{n1}(\tau)(i)=
                    \left\{%
                        \begin{array}{ll}
                            -n, & \hbox{i=1;} \\
                            \tau(i-1), & \hbox{$2 \leq i \leq n.$} \\
                        \end{array}%
                    \right.$
                \item $\varphi_{n2}(\tau)(i)=
                    \left\{%
                        \begin{array}{ll}
                            \tau(i), & \hbox {$1 \leq i \leq n-1$;} \\
                            n, & \hbox{$i=n$.} \\
                        \end{array}%
                    \right.$
            \end{enumerate}
        \end{df}

        Clearly, for $i \in {1,2}, \ \varphi_{ni}$ is a bijection
        (for $U_{n-1}$ onto $U_{ni}$). Therefore, for all
        $\tau_1 \in U_{n1},\ \tau_2 \in U_{n2}$
        there exist unique $\tau_1 ' \in U_{n-1},\ \tau_2 ' \in
        U_{n-1}$, such that $\varphi_{n1}(\tau_1')=\tau_1 , \
        \varphi_{n2}(\tau_2')=\tau_2$. The relations between
        the permutation statistics of $\tau_1, \ \tau_2$, and
        $\tau_1', \ \tau_2'$
        are presented in the following equations:
        \begin{fac}
             \begin{enumerate}
                \item   $neg(\tau_1)=neg(\varphi_{n1}(\tau_1'))=neg(\tau_1')+1$;
                        $neg(\tau_2)=neg(\varphi_{n2}(\tau_2 '))=neg(\tau_2 ')$.
                \item   $sign(\varphi_{n1}(\tau_1 '))=(-1)^n sign(\tau_1')$;
                        $sign(\varphi_{n2}(\tau_2'))=sign(\tau_2')$.
            \end{enumerate}
        \end{fac}
        Now after we finished our preparations we can move on to the
        theorems:

    \medskip

    \subsection{Signed Mahonian}
        \begin{thm}\label{BSF}
                 \[\sum\limits_{\sigma \in B_{n}} sign(\sigma) q^{flag-major_N(\sigma)}=
                    (q;-1)_n [n]_{\pm q^2}!.\]
        \end{thm}

        \noindent{\bf Proof.} We will use the facts mentioned above in \emph{natural linear order} to reduce the summation from $B_n$
        into $U_n$:
        \begin{eqnarray*}
                \sum\limits_{\sigma \in B_n}{sign(\sigma) q^{flag-major_N(\sigma)}}
            &=& \sum\limits_{\tau \pi \in B_n} {sign(\tau \pi) q^{flag-major_N(\tau\pi)}}\\
            &=& \sum\limits_{\tau \in U_n,\ \pi \in S_n}{sign(\pi)sign(\tau) q^{2maj(\pi)+neg(\tau)}}\\
            &=& \sum\limits_{\tau \in U_n} {sign(\tau) q^{neg(\tau)}}\sum\limits_{\pi \in S_n} {sign(\pi) q^{2maj(\pi)}}.
        \end{eqnarray*}
        Now according to Theorem~\ref{GS} we know that:
        \[ \sum_{\pi \in S_n}{sign(\pi)q^{2maj(\pi)}=[1]_{q^2}[2]_{-q^2}[3]_{q^2}[4]_{-q^2}\ldots
        [n]_{(-1)^{n-1}q^2}}=[n]_{\pm q^2}! .\]
        In order to finish the calculation we need to calculate:
        \[a_n=\sum_{\tau \in U_n} {sign(\tau)q^{neg(\tau)}},\]
        we will do so by recursion:
        \begin{eqnarray*}
            a_n &=& \sum\limits_{\tau \in U_n} {sign(\tau) q^{neg(\tau)}}\\
                &=& \sum\limits_{\tau \in U_{n1}} {sign(\tau) q^{neg(\tau)}}+
                    \sum\limits_{\tau \in U_{n2}} {sign(\tau) q^{neg(\tau)}}\\
                &=& \sum\limits_{\tau ' \in U_{n-1}} {sign(\varphi_{n1}(\tau')) q^{neg(\varphi_{n1}(\tau'))}}+
                    \sum\limits_{\tau ' \in U_{n-1}} {sign(\varphi_{n2}(\tau ')) q^{neg(\varphi_{n2}(\tau'))}}\\
                &=& \sum\limits_{\tau ' \in U_{n-1}} {(-1)^n sign(\tau') q^{neg(\tau')+1}}+
                    \sum\limits_{\tau ' \in U_{n-1}} {sign(\tau') q^{neg(\tau')}}\\
                &=&(-1)^nqa_{n-1}+a_{n-1}.
        \end{eqnarray*}
        In the end of the calculation we got the following recurrence equation:
            \[a_n=(1+(-1)^nq)a_{n-1},\ a_1=1-q.\]
        This is easy to solve recurrence equation, which its solution is:
        \begin{eqnarray*}
            a_n=
            \left\{%
                \begin{array}{ll}
                    a_{2m}=(1-q)^m(1+q)^m=(1-q^2)^m, & \hbox{n=2m;} \\
                    a_{2m+1}=(1-q)^{m+1}(1+q)^m=(1-q^2)^m(1-q), & \hbox{n=2m+1.} \\
                \end{array}%
            \right.
        \end{eqnarray*}
        And with the previous calculation we conclude that:
        \begin{eqnarray*}
        \sum_{\sigma \in B_{n}} sign(\sigma)
        q^{flag-major_N(\sigma)}&=&
        \left\{%
            \begin{array}{ll}
                (1-q^2)^m[n]_{\pm q^2}! , & \hbox{n=2m;} \\
                (1-q)(1-q^2)^m[n]_{\pm q^2}!  , & \hbox{n=2m+1.} \\
                \end{array}%
            \right.\\
            &=&(q;-1)_n [n]_{\pm q^2} \quad \bullet
        \end{eqnarray*}
        We move on and calculate the generating function of another \emph{signed Mahonian}
        statistic:
        \begin{thm}\label{signmaj}
            \[ \sum_{\sigma \in B_n} sign(\sigma)q^{nmaj_N(\sigma)}=(q;-q)_n[n]_{\pm q}!.\]
        \end{thm}

        \noindent{\bf Proof.} We will use the same methods as in the
        previous theorem, reducing the problem to $U_n$:

        \begin{eqnarray*}
                \sum\limits_{\sigma \in B_n}sign(\sigma) q^{nmaj_N(\sigma)}&=&
                \sum\limits_{\pi \in S_n, \ \tau \in U_n} sign(\tau\pi) q^{maj(\tau\pi)+nsum(\tau\pi)}\\
            &=& \sum\limits_{\pi \in S_n, \ \tau \in U_n}sign(\pi)sign(\tau)q^{maj(\pi)+nsum(\tau)}\\
            &=& \sum\limits_{\pi \in S_n}sign(\pi)q^{maj(\pi)} \sum\limits_{\tau \in U_n} {sign(\tau)q^{nsum(\tau)}}.
        \end{eqnarray*}
        By Theorem~\ref{GS} we know that $\sum\limits_{\pi \in
        S_n}sign(\pi)q^{maj(\pi)}=[n]_{\pm q}!$, and:
        \begin{eqnarray*}
            a_n &=& \sum\limits_{\tau \in U_n}{sign(\tau)q^{nsum(\tau)}}\\
                &=& \sum\limits_{\tau \in U_{n1}}{sign(\tau)q^{nsum(\tau)}}
                 +  \sum\limits_{\tau \in U_{n2}} {sign(\tau)q^{nsum(\tau)}}\\
                &=& \sum\limits_{\tau' \in U_{n-1}}{sign(\varphi_{n1}(\tau'))q^{nsum(\varphi_{n1}(\tau'))}}
                 +  \sum\limits_{\tau' \in U_{n-1}}{sign(\varphi_{n2}(\tau'))q^{nsum(\varphi_{n2}(\tau'))}}\\
                &=& \sum\limits_{\tau' \in U_{n-1}}{sign(\tau')(-1)^{n}q^{nsum(\tau')+n}}
                 +  \sum\limits_{\tau' \in U_{n-1}}{sign(\tau')q^{nsum(\tau')}}\\
                &=&  (-1)^{n}q^{n}a_{n-1}+a_{n-1}=(1+(-q)^n)a_{n-1}.
        \end{eqnarray*}
        We got a recurrence equation  $a_n=(1+(-1)^nq^n)a_{n-1},\ a_1=1-q,$
        which its solution is:
        \[a_n=\prod\limits_{i=1}^n(1+(-q)^i)=(q;-q)_n.\]
        And after we multiply our results we get:
        \[ \sum_{\sigma \in B_n} sign(\sigma)q^{nmaj_N(\sigma)}=\prod_{i=1}^n(1+(-q)^i)[n]_{\pm q}!=
        (q;-q)_n[n]_{\pm q}! \quad\bullet\]

    \medskip

    \subsection{Bivariate Generating Functions}
        We can use the method of reducing the summation form $B_n$ into $U_n$, for the purpose of
        generalizing and proving known theorems, we will start with
        generalizing Theorem~\ref{lenmaj}:

        \begin{thm}
            \[\sum_{\sigma \in B_n} q^{nmaj_N(\sigma)}t^{neg(\sigma)}=\sum_{\sigma \in B_n}
            q^{\ell(\sigma)}t^{neg(\sigma)}.\]
        \end{thm}

        \noindent{\bf Proof.}
        \begin{eqnarray*}
            \sum\limits_{\sigma \in B_n}q^{nmaj_N(\sigma)}t^{neg(\sigma)}&=&
                \sum\limits_{\pi \in S_n,\ \tau \in U_n} q^{maj(\tau \pi)+nsum(\tau\pi)}t^{neg(\tau \pi)}\\
            &=& \sum\limits_{\pi \in S_n,\ \tau \in U_n} q^{maj(\pi)+nsum(\tau)}t^{neg(\tau)}\\
            &=& \sum\limits_{\tau \in U_n} q^{nsum(\tau)}t^{neg(\tau)}\sum\limits_{\pi \in S_n} q^{maj(\pi)}\\
            &=& \sum\limits_{\tau \in U_n} q^{nsum(\tau)}t^{neg(\tau)}\sum\limits_{\pi \in S_n}q^{inv(\pi)}\\
            &=& \sum\limits_{\pi \in S_n,\ \tau \in U_n}q^{inv(\tau\pi)+nsum(\tau\pi)}t^{neg(\tau\pi)}\\
            &=& \sum\limits_{\sigma \in B_n}q^{\ell(\sigma)}t^{neg(\sigma)} \quad\bullet
        \end{eqnarray*}
        We can also calculate some interesting \emph{Mahonian-Mahonian} generating function:
        \begin{thm}\label{fnmaj}
            \[\sum_{\sigma \in B_n} q^{flag-major_N(\sigma)}t^{nmaj_N(\sigma)}=\prod_{i=1}^n(1+qt^i)[n]_{q^2t}!.\]
        \end{thm}

        \noindent{\bf Proof.} We will prove this theorem by reducing the
        problem to $U_n$:

        \begin{eqnarray*}
            \sum_{\sigma \in B_n} q^{flag-major_N(\sigma)}t^{nmaj_N(\sigma)}&=&
                \sum_{\pi \in S_n,\ \tau \in U_n}{q^{2maj(\pi)+neg(\tau)}t^{maj(\pi)+nsum(\tau)}}\\
            &=& \sum_{\tau \in U_n}{q^{neg(\tau)}t^{nsum(\tau)}}\sum_{\pi \in S_n}{q^{2maj(\pi)}t^{maj(\pi)}}\\
            &=& \sum_{\tau \in U_n}{q^{neg(\tau)}t^{nsum(\tau)}}\sum_{\pi \in S_n}{(q^2t)^{maj(\pi)}}.
        \end{eqnarray*}
        We know according to Theorem~\ref{MM} that:
        $\sum\limits_{\pi \in S_n}{(q^2t)^{maj(\pi)}}=[n]_{q^2t}!$, and by calculation we get:

        \begin{eqnarray*}
            a_n=\sum_{\tau \in U_n} {q^{neg(\tau)}t^{nsum(\tau)}}&=&
                \sum_{\tau \in U_{n1}} {q^{neg(\tau)}t^{nsum(\tau)}}
             +  \sum_{\tau \in U_{n2}}q^{neg(\tau)}t^{nsum(\tau)}\\
            &=& \sum_{\tau' \in U_{n-1}}
                {q^{neg(\varphi_{n1}(\tau'))}t^{nsum(\varphi_{n1}(\tau'))}}\\
            &+&  \sum_{\tau' \in U_{n-1}}q^{neg(\varphi_{n2}(\tau'))}t^{nsum(\varphi_{n2}(\tau'))}\\
            &=& \sum_{\tau' \in U_{n-1}} {q^{neg(\tau')+1}t^{nsum(\tau')+n}}
             +  \sum_{\tau' \in U_{n-1}}q^{neg(\tau')}t^{nsum(\tau')}\\
            &=& qt^{n}a_{n-1}+a_{n-1}=(1+qt^n)a_{n-1}.
        \end{eqnarray*}
        We got the recurrence equation: $a_n=(1+qt^n)a_{n-1},\ a_1=1+qt$, and
        the solution to this equation is:
        $a_n=\prod\limits_{i=1}^n(1+qt^i)$, and therefore; the general solution is:
        \[\sum_{\sigma \in B_n} q^{flag-major_N(\sigma)}t^{nmaj_N(\sigma)}=[n]_{q^2t}!\prod_{i=1}^n(1+qt^i) \quad\bullet\]

        \begin{nt} Notice that  if we put $t=1$ in the previous Theorem~\ref{fnmaj}, we
            get Theorem~\ref{fmajlen} and the equation: $[n]_{q^2}(1+q)^n=\prod\limits_{i=1}^n [2i]_q!$.
        \end{nt}
        We can also calculate the generating function of \emph{length} and \emph{flag major index} by using a similar method:
        \begin{thm}\label{fmajlenN}
            \[\sum_{\sigma \in B_n} q^{flag-major_N(\sigma)}t^{\ell(\sigma)}=
                A_n(q^2,t)\prod_{i=1}^n (1+qt^i),\] where
                $A_n(q,t)=\sum\limits_{\pi \in S_n} q^{maj(\pi)}t^{\ell(\pi)}=
                \sum\limits_{\pi \in S_n} q^{maj(\pi)}t^{inv(\pi)}.$
        \end{thm}

        \noindent{\bf Proof.} We will start with reducing $B_n$ to
        $U_n$:
        \begin{eqnarray*}
            \sum\limits_{\sigma \in B_n}q^{flag-major(\sigma)}t^{\ell(\sigma)}&=&
                \sum\limits_{\tau \in U_n, \ \pi \in S_n} q^{2maj(\pi)+neg(\tau)}t^{inv(\pi)+nsum(\tau)}\\
            &=& \sum\limits_{\tau \in U_n} q^{neg(\tau)}t^{nsum(\tau)} \sum\limits_{\pi \in S_n} q^{2maj(\pi)}t^{inv(\pi)}\\
            &=& A_n(q^2,t)\sum\limits_{\tau \in U_n} q^{neg(\tau)}t^{nsum(\tau)}.
        \end{eqnarray*}
        Now we need to calculate $\sum\limits_{\tau \in U_n}
        q^{neg(\tau)}t^{nsum(\tau)}$, we calculate this exact
        equation in the proof of the previous Theorem~\ref{fnmaj},
        and got the recurrence equation:
        $a_n=(1+qt^n)a_{n-1},\ a_1=1+qt$, which its solution is:
        $a_n=\prod\limits_{i=1}^n (1+qt^i)$, therefore; the
        final solution is:
        \[\sum_{\sigma \in B_n} q^{flag-major_N (\sigma)}t^{\ell(\sigma)}=
                \prod_{i=1}^n (1+qt^i)A_n(q^2,t) \quad \bullet \]

    \medskip

    \subsection{The Signed Mahonian Over $D_n$}
            In this subsection we are going to calculate \emph{signed Mahonian} on
            $D_n$ the subgroup of $B_n$ (see definition~\ref{wsubgroups}).
            We open this subsection with few definitions that we are going to use in our proofs:
            \begin{df} Define the following subset of $D_n$:
                \begin{eqnarray*}
                    UD_n:=U_n \cap D_n=\{\tau \in D_n \ |\ \tau(1)<\tau(2)<\ldots<\tau(n)\}.
                \end{eqnarray*}
            \end{df}

            \begin{rem} Let be $\tau \in UD_n$,
                if $\tau(n)=n$, then exists unique $\tau' \in UD_{n-1}$, where $\varphi_{n1}(\tau')=\tau$,
                and if $\tau(1)=-n$ then exists unique $\tau' \in \overline{UD_{n-1}}:=U_{n-1}-UD_{n-1},$
                where $\varphi_{n2}(\tau')=\tau$.
            \end{rem}
            Notice that for every $\sigma \in D_n$ there exists a unique decomposition:
            $\sigma=\tau\pi, \ \pi \in S_n, \ \tau \in UD_n.$ Now we use the last
            remark, and the constructions in the beginning of this subsection to prove the following theorem:

            \begin{thm} \label{DSF}
                \begin{eqnarray*}
                    \sum\limits_{\sigma \in D_n} sign(\sigma)q^{flag-major_N(\sigma)}=
                    (1-q^2)^{\lfloor\frac{n}{2}\rfloor}[n]_{\pm q^2}!.
                \end{eqnarray*}
            \end{thm}

            \noindent{\bf Proof.}
            \begin{eqnarray*}
                    \sum\limits_{\sigma \in D_n} sign(\sigma)q^{flag-major_N(\sigma)}&=&
                    \sum\limits_{\pi \in S_n ,\ \tau\in UD_n} sign(\tau)sign(\pi)q^{2maj(\pi)+neg(\tau)}\\
                &=& \sum\limits_{\pi \in S_n} sign(\pi)q^{2maj(\pi)} \sum\limits_{\tau \in UD_n} sign(\tau)q^{neg(\tau)}\\
                &=& [n]_{\pm q^2}! \sum\limits_{\tau \in UD_n}sign(\tau)q^{neg(\tau)}.
            \end{eqnarray*}
            Now we need to calculate $\sum\limits_{\tau\in
            UD_n}sign(\tau)q^{neg(\tau)}$, we will do so by recursion to
            $UD_{n-1}$, before doing so, we notice the fact that:
            \[  \sum_{\tau \in U_n}sign(\tau)q^{neg(\tau)}=\sum_{\tau \in UD_n}sign(\tau)q^{neg(\tau)}+
            \sum_{\tau \in \overline{UD_n}}sign(\tau)q^{neg(\tau)},\]
            therefore;

            \begin{eqnarray*}
                \sum\limits_{\tau \in \overline{UD_n}}sign(\tau)q^{neg(\tau)}&=&
                    \sum\limits_{\tau \in U_n}sign(\tau)q^{neg(\tau)}
                -   \sum_{\tau \in UD_n}sign(\tau)q^{neg(\tau)}\\
                &=&  \left\{%
                        \begin{array}{ll}
                            (1-q^2)^m-\sum\limits_{\tau \in UD_{2m}}sign(\tau)q^{neg(\tau)}, & \hbox{n=2m} \\
                            (1-q)(1-q^2)^m -\sum\limits_{\tau \in UD_{2m+1}}sign(\tau)q^{neg(\tau)}  , & \hbox{n=2m+1.}
                        \end{array}%
                    \right.
            \end{eqnarray*}
            Now we are ready to do the calculation:
            \begin{eqnarray*}
                a_n &=&\sum_{\tau \in UD_n} sign(\tau)q^{neg(\tau)}\\
                &=& \sum_{\tau \in UD_n , \ \tau(n)=n} sign(\tau)q^{neg(\tau)}
                 +  \sum_{\tau \in UD_n , \ \tau(1)=-n} sign(\tau)q^{neg(\tau)}\\
                &=& \sum_{\tau' \in UD_{n-1}} sign(\varphi_{n2}(\tau'))q^{neg(\varphi_{n2}(\tau'))}
                 +  \sum_{\tau' \in \overline{UD_{n-1}}} sign(\varphi_{n1}(\tau'))q^{neg(\varphi_{n1}(\tau'))}\\
                &=& \sum_{\tau' \in UD_{n-1}} sign(\tau')q^{neg(\tau')}
                 +  \sum_{\tau' \in \overline{UD_{n-1}}}(-1)^{n}sign(\tau')q^{neg(\tau')+1}\\
                &=&  a_{n-1}+(-1)^nq[(1-q)^{\lceil \frac{n-1}{2} \rceil} (1+q)^{\lfloor \frac{n-1}{2}
                    \rfloor}-a_{n-1}].
            \end{eqnarray*}
            We got the recurrence equations:
            \begin{eqnarray*}
                    a_n=
                    \left\{%
                        \begin{array}{ll}
                            (1-q)a_{2m-1}+q(1-q^2)^{m-1}(1-q), & \hbox{n=2m;} \\
                            (1+q)a_{2m}-q(1-q^2)^{m}, & \hbox{n=2m+1.} \\
                        \end{array}%
                    \right.
                , a_1=1.
            \end{eqnarray*}
            The solutions to these equations are:
            \begin{eqnarray*}
                a_{2m+1}&=&(1+q)a_{2m}-q(1-q^2)^m\\
                &=& (1+q)[(1-q)a_{2m-1}+q(1-q^2)^{m-1}(1-q)]-q(1-q^2)^m\\
                &=& (1-q^2)a_{2m-1}+q(1-q^2)^m-q(1-q^2)^m=(1-q^2)a_{2m-1}\\
                &=& (1-q^2)^m.\\ \\
                a_{2m}&=&(1-q)(1-q^2)^{m-1}+q(1-q)(1-q^2)^{m-1}\\
                &=& (1-q)(1-q^2)^{m-1}(1+q)\\
                &=&(1-q^2)^m.
            \end{eqnarray*}
            In the end of the equation we got the solution:
            \begin{eqnarray*}
                a_n=(1-q^2)^{\lfloor\frac{n}{2}\rfloor}=
                \left\{%
                    \begin{array}{ll}
                        (1-q^2)^m, & \hbox{n=2m;} \\
                        (1-q^2)^m, & \hbox{n=2m+1.} \\
                    \end{array}%
                \right.
            \end{eqnarray*}
            Now we can write the final solution:
            \begin{eqnarray*}
                \sum\limits_{\sigma \in D_n} sign(\sigma)q^{flag-major_N(\sigma)}=
                (1-q^2)^{\lfloor\frac{n}{2}\rfloor}[n]_{\pm q^2}!
                \quad \bullet
            \end{eqnarray*}

        \begin{cor}\label{DSFC} From the last Theorem~\ref{DSF}, and from Theorem~\ref{BSF}
            we can conclude that for even $n$:
            \[ \sum\limits_{\sigma \in D_n} sign(\sigma)q^{flag-major_N(\sigma)}=
            \sum\limits_{\sigma \in B_n} sign(\sigma)q^{flag-major_N(\sigma)}.\]
        \end{cor}

\section{Colored Permutation Groups}
    In this section we present a method, which takes known equalities from $S_n$
    and generalize them to $G_{r,n}$. The main idea of this method is to rely on the fact that our
    alphabet:
    \[\{1^{[r-1]},\ldots,1,2^{[r-1]},\ldots,2,\ldots,n^{[r-1]},\ldots,n\}\]
    is organized according to the \emph{friends order}:
    \[F: \ 1^{[r-1]}<\ldots<1<2^{[r-1]}<\ldots<2<\ldots<n^{[r-1]}<\ldots<n.\]
    In the \emph{friends order} we can notice, that the colors have no effect on many statistics.

    \subsection{Preliminaries}
    We begin with the definition of the \emph{$Color^n_r$} groups:
    \begin{df}\label{COLORN} Let $r,n$ be nonnegative integers.
        Define the $Color^n_r$ group to be:
        \[Color^n_r:=\{(z,id) \in G_{r,n} \}\cong \mathbf{Z}^n_r.\]
    \end{df}
    We define the following subsets of $Color^n_r$:
    \begin{df}  Let $i \in [1,r]$. Define:
        \[Color_{r,i}^n:=\{(z,id) \in Color^n_r \ | \ z_n=i\}.\]
    \end{df}

    \begin{rem}
        The group $Color^n_r$ is a disjoint union:
        \[Color^n_r=\biguplus\limits_{i=0}^{r-1} Color_{r,i}^n,\]
        where $\biguplus$ donates for disjoint union.
    \end{rem}
    \label{taupi}

    Every $(z,\pi) \in G_{r,n}$ has two a unique compositions as
    products $(z,\pi)=(z,id)(0,\pi)=(0,\pi')(z',\pi'')$, where here
    $0:=(0,\ldots,0)$.

    \begin{lem}\label{tpc} Let $i \in [1,k]$ and let $ostat_i: \ S_n \to \mathbf{Z}, \ cstat_i: \ \mathbf{Z}^n_r \to \mathbf{Z}$
        statistics on $S_n$ and $\mathbf{Z}^n_r$ respectively; define $stat_i: G_{r,n} \to
        \mathbf{Z}$ by $stat_i(z,\pi):=ostat_i(\pi)+cstat(z)$. Then:
        \begin{eqnarray*}
            \sum_{(z,\pi) \in G_{r,n}} q_1^{stat_1(z,\pi)}
            \ldots q_k^{stat_k(z,\pi)}=
            \sum_{(z,id) \in Color^n_r} q_1^{cstat_1(z)}
            \ldots q_k^{cstat_k(z)}\sum_{\pi \in S_n} q_1^{ostat_1(\pi)}
            \ldots q_k^{ostat_k(\pi)}.
        \end{eqnarray*}
    \end{lem}
    \noindent{\bf Proof.} By using the fact that every
    $(z,\pi)=(z,id)(0,\pi),$
    \begin{eqnarray*}
        \sum\limits_{(z,\pi) \in G_{r,n}} q_1^{stat_1(z,\pi)}
                \ldots q_k^{stat_k(z,\pi)}&=&
        \sum\limits_{(z,\pi) \in G_{r,n}} q_1^{ostat_1(z,\pi)+cstat_1(z,\pi)}
                \ldots q_k^{ostat_k(z,\pi)+cstat_k(z,\pi)}\\
        &=&     \sum\limits_{ (z,id) \in Color^n_r}\sum\limits_{(0,\pi) \in G_{r,n}} q_1^{ostat_1(\pi)+cstat_1(z)}
                \ldots q_k^{ostat_k(\pi)+cstat_k(z)}\\
        &=&     \sum\limits_{(z,id) \in Color^n_r}
                q_1^{cstat_1(z)} \ldots q_k^{cstat_k(z)}
                 \sum\limits_{(0,\pi) \in G_{r,n}}
                q_1^{ostat_1(\pi)} \ldots q_k^{ostat_k(\pi)}\\
        &=&     \sum\limits_{(z,id) \in Color^n_r}
                q_1^{cstat_1(z)} \ldots q_k^{cstat_k(z)}
                \sum\limits_{\pi \in S_n}
                q_1^{ostat_1(\pi)} \ldots q_k^{ostat_k(\pi)} \bullet
    \end{eqnarray*}
    We now define a set of bijections $Color^{n-1}_r \to
    Color_{r,j}^n, \ j \in [0,r-1]$:
    \begin{df} Let $j \in [0,r-1]$. Define the bijections\\ $\xi^j: Color^{n-1}_r \to
        Color_{r,j}^n$ by:
        \begin{eqnarray*}
             \xi^j((z,id)):=(z',id), \mbox{ where }
                \left\{%
                    \begin{array}{ll}
                        z'_i=z_i, &\hbox{$1 \leq i \leq n-1$;} \\
                        z'_i=j,  &\hbox{$i=n$.}
                    \end{array}%
                \right. \ \ ( \forall (z,id)\in Color^{n-1}_r).
        \end{eqnarray*}
    \end{df}
    The relations between $(z,id)$, and $\xi^j((z',id))$ are presented
    in the following equation that can be easily concluded from the
    definition of $\xi^j$:

    \begin{enumerate}
        \item $ncsum(\xi^j((z,id)))=
                \left\{%
                    \begin{array}{ll}
                        ncsum((z,id))+n-1, & \hbox{$1 \leq j \leq r-1;$} \\
                        ncsum((z,id)), & \hbox{j=0.} \\
                    \end{array}%
                \right.$\\
                Recall that $ncsum((z,\pi)):=\sum\limits_{z_i \neq 0}
                    (\pi(i)-1).$
        \item $csum(\xi^j((z,id)))=csum((z,id))+j.$
    \end{enumerate}

    \subsection{Friends Color-Signed Mahonian}
        In this subsection we calculate the bivariate generating
        function of the \emph{flag major index} and \emph{friends color sign}
        statistics. We open this subsection with the definition of
        the \emph{color sign} statistic:
        \begin{df}\label{csign} Let $r,n$ be nonnegative integers. Define the
            \emph{color sign} of $g \in G_{r,n}$ to be:
            \[csign(g):=(-1)^{\ell(g)}.\]
        \end{df}

        \begin{nt} Notice that in the general case of \emph{r-colors} the \emph{color-sign} statistic
            is not necessarily a character. For example in the case of $n=3,\ r=3$, we
            take the following permutations:
             \begin{eqnarray*}
                &&  g_1=((1,0,0),[1,2,3]), \ g_2=((2,0,0),[1,2,3])
                    \mbox{ therefore; } g_1g_2=((0,0,0),[1,2,3]),
                    \mbox{ and }\\
                &&  csign(g_1)=-1, \ csign(g_2)=1,
                \mbox{ but } csign(g_1)csign(g_2)=-1\neq
                csign(g_1g_2)=1.
             \end{eqnarray*}
        \end{nt}
        \begin{df}\label{fcsign} Let $g \in G_{r,n}$. Define the  \emph{friends color-sign} of $g \in G_{r,n}$ to be:
            \[csign_{F}(g):=(-1)^{inv_{F}(g)+csum(g)}.\]
        \end{df}
        \begin{pro} Let $g=(z,\pi) \in G_{r,n}$ then:
            \[csign_F(g)=csign_F((0,\pi))csign_F((z,id)).\]
        \end{pro}
        \noindent{\bf Proof.} According to Definition~\ref{fcsign}
        and Lemma~\ref{tpc} we conclude:
        \begin{eqnarray*}
            csign_F((z,\pi))&=&(-1)^{inv_F((z,\pi))+csum((z,\pi))}\\
            &=&(-1)^{inv_F((0,\pi))+csum((z,id))}\\
            &=&(-1)^{inv_F((0,\pi))+csum((0,\pi))}(-1)^{inv_F((z,id))+csum((z,id))}\\
            &=& csign_F((0,\pi))csign_F((z,id)) \quad \bullet
        \end{eqnarray*}

        Now by using the above facts, we can calculate the generating function of
        \emph{friends color sign} and \emph{flag major index}:
        \begin{thm}\label{csignfmaj}
            \[\sum_{g \in G_{r,n}} csign_F(g)q^{flag-major_F(g)}= [r]_{-q}^n [n]_{\pm rq}!.\]
        \end{thm}

        \noindent{\bf Proof.} According to Lemma~\ref{tpc} and Theorem~\ref{GS} we get:
        \begin{eqnarray*}
            \sum\limits_{g \in G_{r,n}}
                csign_F(g)q^{flag-major_F(g)}&=&
            \sum\limits_{(z,id) \in Color^n_r} csign_F((z,id))q^{csum((z,id))} \sum\limits_{\pi \in S_n} sign(\pi)q^{r\cdot maj(\pi)} \\
            &=& \sum\limits_{(z,id) \in Color^n_r}
                    csign_F((z,id)) q^{csum((z,id))}[n]_{\pm q^r}!.
        \end{eqnarray*}
        Now we need to calculate: $a_n=\sum\limits_{(z,id) \in Color^n_r }
        csign_F((z,id)) q^{csum((z,id))}$.
        We will do so by using the $\xi^j$ bijections to find recurrence equation to $a_n$:
        \begin{eqnarray*}
            a_n&=&
                \sum\limits_{(z,id) \in Color^n_r}
                    csign_F((z,id))q^{csum((z,id))}\\
            &=& \sum\limits_{i=0}^{r-1}\sum\limits_{(z,id) \in
                    Color_{r,i}^n} csign_F((z,id))q^{csum((z,id))}\\
            &=& \sum\limits_{i=0}^{r-1}\sum\limits_{(z,id) \in
                    Color_r^{n-1}} csign_F(\xi^i((z,id)))q^{csum(\xi^j((z,id)))}\\
            &=&      \sum\limits_{i=0}^{r-1}\sum\limits_{(z,id) \in
                        Color_r^{n-1}} (-1)^{i} csign_F((z,id))q^{csum((z,id))+i}\\
            &=& \sum\limits_{i=0}^{r-1} (-q)^ia_{n-1}\\
            &=& [r]_{-q}a_{n-1}.
        \end{eqnarray*}
        In the end we got the recurrence equation:
            \[a_n=[r]_{-q}a_{n-1}, \ a_1=[r]_{-q},\]
        which its solution is:
        \[a_n=[r]_{-q}^n.\]
        Our final solution is:
        \[\sum_{g \in G_{r,n}} csign_F(g)q^{flag-major_F(g)}=
            [r]_{-q}^n[n]_{\pm rq}! \quad \bullet\]

    \begin{nt} Notice that if we put $r=1$ in Theorem~\ref{csignfmaj}, we get
        Theorem~\ref{GS}.
    \end{nt}

    \medskip

    \subsection{Colored Excedance and Denert's Statistics}
        In this subsection we give a generalized version of the \emph{excedance number} and \emph{Denert's statistic}
        for the \emph{colored permutation groups}, by using the \emph{r-color excedance number} (denoted $exc_r$) and the
        \emph{r-color Denert's statistic} (denoted $den_r$), similar type results appear also in \cite{ST}.

        \begin{df}\label{excr} Let $r,n$ be a nonnegative integers and $g=(z,\pi) \in G_{r,n}$.
         Define the \emph{r-color excedance number} of $g$ to be:
            \[exc_r(g):=exc(\pi)+csum(g).\]
        \end{df}
        We prove that the \emph{r-color excedance number} is \emph{Eulerian}:
        \begin{thm} Let $r,n$ be nonnegative integers. Then:
            \[\sum_{g \in G_{r,n}} q^{exc_r(g)}=\sum_{g \in
            G_{r,n}} q^{ldes_F(g)}.\]
        \end{thm}
        \noindent{\bf Proof.} According to Lemma~\ref{tpc} and Theorem~\ref{desexcA} we get:
        \begin{eqnarray*}
            \sum\limits_{g \in G_{r,n}} q^{exc_r(g)}&=&
                \sum\limits_{g \in G_{r,n}} q^{exc(\pi)+csum(g)}\\
            &=& \sum\limits_{(z,id) \in Color^n_r} q^{csum((z,id))}
                \sum\limits_{\pi \in S_n} q^{exc(\pi)}\\
            &=& \sum\limits_{(z,id) \in Color^n_r} q^{csum((z,id))}
                \sum\limits_{\pi \in S_n} q^{des(\pi)}\\
            &=& \sum\limits_{(z,id) \in Color^n_r} \sum\limits_{\pi \in S_n} q^{des(\pi)+csum((z,id))}\\
            &=& \sum_{g \in G_{r,n}} q^{ldes_F(\sigma)} \qquad \bullet
        \end{eqnarray*}

        \begin{df}\label{denr} Let $r,n$ be nonnegative integers.
            Define for $g \in G_{r,n}$ the \emph{r-color Denert's statistic} to be:
            \[den_r(g)=r\cdot den(\pi)+csum(g).\]
        \end{df}
        We prove that the \emph{r-color Denert's statistic} is
        equidistributed with the \emph{flag major index} over $G_{r,n}$:
        \begin{thm} Let $r,n$ be a nonnegative integers. Then:
            \[\sum_{g \in G_{r,n}} q^{den_r(g)}=\sum_{g \in
            G_{r,n}} q^{flag-major_F(g)}.\]
        \end{thm}
        \noindent{\bf Proof.} Using Lemma~\ref{tpc}, and
        Theorem~\ref{FZ} we get:
        \begin{eqnarray*}
            \sum\limits_{g \in G_{r,n}} q^{den_r(g)}&=&
                \sum\limits_{g \in G_{r,n}} q^{r\cdot den(\pi)+csum(g)}\\
            &=& \sum\limits_{(z,id) \in Color^n_r} q^{csum((z,id)))}
                \sum\limits_{\pi \in S_n} q^{r\cdot den(\pi)}\\
            &=& \sum\limits_{(z,id) \in Color^n_r} q^{csum((z,id)))}
                \sum\limits_{\pi \in S_n} q^{r\cdot maj(\pi)}\\
            &=& \sum\limits_{(z,id) \in Color^n_r}\sum\limits_{\pi \in S_n}  q^{r\cdot maj(\pi)+csum((z,id))}\\
            &=& \sum_{g \in G_{r,n}} q^{flag-major_F(g)} \qquad \bullet
        \end{eqnarray*}
        Now we prove that the pair of statistics $(den_r,exc_r)$
        is equidistributed with $(flag-major_F,ldes_F)$ and therefore,
        the pair of statistics $(den_r,exc_r)$ generalizes the Foata-Zeilberger Theorem~\ref{FZ}.
        \begin{thm}\label{FZR} Let $r,n$ be a nonnegative integers. Then:
            \[\sum_{g \in G_{r,n}} q^{den_r(g)}t^{exc_r(g)}=\sum_{g \in
            G_{r,n}} q^{flag-major_F(g)}t^{ldes_F(g)}.\]
        \end{thm}
        \noindent{\bf Proof.} Using Lemma~\ref{tpc} and Theorem~\ref{FZ} we can
        conclude:
        \begin{eqnarray*}
            \sum\limits_{g \in G_{r,n}} q^{den_r(g)}t^{exc_r(g)}&=&
                \sum\limits_{g \in G_{r,n}} q^{r\cdot den(\pi)+csum(g)}t^{exc(g)+csum(g)}\\
            &=& \sum\limits_{(z,id) \in Color^n_r} q^{csum((z,id))}t^{csum((z,id))}
                  \sum\limits_{\pi \in S_n} q^{r\cdot den(\pi)}t^{exc(\pi)}\\
            &=& \sum\limits_{(z,id) \in Color^n_r} q^{csum((z,id))}t^{csum((z,id))}
                    \sum\limits_{\pi \in S_n} q^{r\cdot maj(\pi)}t^{des(\pi)}\\
            &=& \sum\limits_{g \in G_{r,n}}
                    q^{r\cdot maj(\pi)+csum(g)}t^{des(g)+csum(g)}\\
            &=& \sum\limits_{g \in G_{r,n}} q^{flag-major_F(g)}t^{ldes_F(g)}\quad \bullet
        \end{eqnarray*}

    \medskip

    \subsection{Flag-Excedance and Flag-Denert's Statistic of Type B}
        In this subsection we present the \emph{flag-Denert's} statistic (denoted $fden$) and
        the \emph{flag-excedance} (denoted $fexc$) statistic. We prove that the pair
        of statistics $(fden,fexc)$ are equidistributed with $(flag-major,fdes)$
        over $B_n$ and, therefore, the \emph{flag-Denert} and \emph{flag-excedance} statistics gives
        a type B generalization to the Foata-Zeilberger Theorem~\ref{FZ}.

        \begin{df} Define the type b \emph{excedance} number of $\sigma \in B_n$ to be:
            \[exc_B(\sigma):=|\{\ 1 \leq i \leq n \ |  \ i<
            |\sigma(i)| \ \}|.\]
        \end{df}

        \begin{df}\label{fexc} Define the \emph{flag-excedance} of $\sigma \in B_n$ to be:
            \[fexc(\sigma):=2exc_B(\sigma)+\varepsilon(\sigma).\]
        \end{df}
        \begin{df} Let $n$ be a nonnegative integer. Define the following
            subset of $B_n$:
                \[Color_2^n:=\{\sigma \in B_n | \ \sigma(i)=\pm i, \ \forall i \in [1,n] \}.\]
            This is another form of definition~\ref{COLORN} for
            $r=2$.
        \end{df}

        We prove that the \emph{flag-excedance} statistics is
        Eulerian.
        \begin{rem}\label{remfexc} In the following theorems we use the fact that
            if $\sigma=\pi\tau$ and where $\pi=[|\sigma(1)|,\ldots,|\sigma(n)|\ ]$ and $\tau \in
            Color_2^n$, then $\varepsilon(\sigma)=\varepsilon(\tau).$
        \end{rem}

        \begin{thm}
            \[\sum_{\sigma \in B_n} q^{fexc(\sigma)}=
            \sum_{\sigma \in B_n} q^{fdes_F(\sigma)}.\]
        \end{thm}

        \noindent{\bf Proof.} By using Lemma~\ref{tpc} and Remark~\ref{remfexc}
        we can conclude that:
        \begin{eqnarray*}
            \sum\limits_{\sigma \in B_n} q^{fexc(\sigma)}&=&
                \sum\limits_{\tau \in Color^n_2}\sum\limits_{ \pi \in S_n}
                    q^{2exc_B(\pi\tau)+\varepsilon(\pi\tau)}\\
            &=& \sum\limits_{\tau \in Color^n_2} q^{\varepsilon(\tau)}\sum\limits_{\pi \in S_n}
                    q^{2exc(\pi)}\\
            &=& \sum\limits_{\tau \in Color^n_2} q^{\varepsilon(\tau)}\sum\limits_{\pi \in S_n}
                    q^{2des(\pi)}\\
            &=& \sum\limits_{\sigma \in B_n} q^{fdes_F(\sigma)}
                \qquad \bullet
       \end{eqnarray*}
       We define the type B Denert's statistic (denoted $den_B$):
       \begin{df} Let $\sigma \in B_n$. Define the type B Denert's statistic to be:
            \begin{eqnarray*}
                den_B(\sigma)&=&|\{1 \leq l<k \leq n \ | \ |\sigma(k)|<|\sigma(l)|<k\}|\\
                &+&|\{1 \leq l<k \leq n \ | \ |\sigma(l)|<k<|\sigma(k)|\}|\\
                &+&|\{1 \leq l<k \leq n \ | \ k<|\sigma(k)|<|\sigma(l)|\}|.
            \end{eqnarray*}
       \end{df}
       \begin{rem}\label{remden} According to the definition of $den_B$ we can
           see that: \[den_B(\sigma)=den_B(\tau\pi)=den_B(\pi), \ \forall \sigma \in B_n, \
            \tau \in Color_2^n, \ \pi \in S_n.\]
       \end{rem}
       We define the \emph{flag-Denert's} statistic (denoted $fden_B$),
       and prove that it is equidistributed with the
       \emph{flag major index} over the \emph{signed permutations
       groups}:

       \begin{df}\label{fden1} Let $\sigma \in B_n$. Define the \emph{flag-Denert's} statistic to be:
            \[ fden(\sigma):=2den_B(\sigma)+neg(\sigma).\]
       \end{df}

       \begin{thm} \[\sum_{\sigma \in B_n} q^{fden(\sigma)}=\sum_{\sigma \in B_n}
            q^{flag-major_F(\sigma)}.\]
       \end{thm}
       \noindent{\bf Proof.} We use the Definition~\ref{fden1}, Lemma~\ref{tpc}, and Remark~\ref{remden}
       and conclude the following equality:
       \begin{eqnarray*}
            \sum\limits_{\sigma \in B_n} q^{fden(\sigma)}&=&
                \sum\limits_{\sigma \in B_n}
                    q^{2den_B(\sigma)+neg(\sigma)}\\
            &=& \sum\limits_{\tau \in Color^n_2} q^{neg(\tau)}
                    \sum\limits_{\pi \in S_n}q^{2den(\pi)}\\
            &=& \sum\limits_{\tau \in Color^n_2} q^{neg(\tau)} \sum\limits_{\pi \in S_n} q^{2maj(\pi)}\\
            &=&  \sum\limits_{\sigma \in B_n} q^{flag-major_F(\sigma)} \qquad \bullet
       \end{eqnarray*}
     We prove that the pair of statistics
     (\emph{fden},\emph{fexc}) is equidistributed with
     (\emph{flag-major},\emph{fdes}).
     \begin{thm}\label{fden}
     \[\sum_{\sigma \in B_n} q^{fden(\sigma)}t^{fexc(\sigma)}=\sum_{\sigma \in B_n}
        q^{flag-major_F(\sigma)}t^{fdes_F(\sigma)}.\]
     \end{thm}
     \noindent{\bf Proof.} We use the Definitions~\ref{fexc},~\ref{fden1}, Lemma~\ref{tpc},
     and Theorem~\ref{FZ} and conclude the following equality:
     \begin{eqnarray*}
            \sum\limits_{\sigma \in B_n}
                q^{fden(\sigma)}t^{fexc(\sigma)}&=&
            \sum\limits_{\sigma \in B_n} q^{2den_B(\sigma)+neg(\sigma)}t^{2exc_B(\sigma)+\varepsilon(\sigma)}\\
        &=& \sum\limits_{\tau \in Color^n_2}q^{neg(\tau)}t^{\varepsilon(\tau)}
                \sum\limits_{\pi \in S_n} q^{2den(\pi)}t^{2exc(\pi)}\\
        &=& \sum\limits_{\tau \in Color^n_2}q^{neg(\tau)}t^{\varepsilon(\tau)}
                \sum\limits_{\pi \in S_n} q^{2maj(\pi)}t^{2des(\pi)}\\
        &=& \sum\limits_{\sigma \in B_n } q^{flag-major_F(\sigma)}t^{fdes_F(\sigma)}
            \qquad\bullet
   \end{eqnarray*}

\section{The Flag-Inversion Statistic}
    In this section we define a new statistic, \emph{flag-inversion} (denoted $finv$).
    We show that this new statistic is also equidistributed with the \emph{flag major
    index}, and therefore it is equidistributed with the \emph{length} statistic
    for $r=1,2$. Using \emph{finv} we calculate interesting generating
    functions and present a generalization of Foata and Sch\"{u}tzenberger's Theorem~\ref{FS1} to the groups
    $G_{r,n}$. We finish this section by giving an algebraic interpretation
    of \emph{finv}.

    \begin{rem} Throughout this section we assume that our
        alphabet is:
            \[\{1,\ldots,1^{[r-1]},2,\ldots,2^{[r-1]},\ldots,\ldots,n,\ldots,n^{[r-1]}\},\]
        and it is ordered by the \emph{friends order}:
        \[F: 1<\ldots<1^{[r-1]}<2<\ldots<2^{[r-1]}<\ldots<n<\ldots<n^{[r-1]}.\]
        For convenience we define $i^{[0]}:=i,\ \forall i \in [1,n].$
    \end{rem}

    \medskip
    \subsection{A New Mahonian Statistic}

    \begin{df}\label{finv} Let $g \in G_{r,n}$. Define the
        \emph{flag-inversion} statistic of $g$ to be
            \[finv(g):=r \cdot inv(g)+csum(g).\]
    \end{df}

    \begin{rem} For $\sigma \in B_n$ this can be written as:
        \[finv(\sigma):=2\cdot inv(\sigma)+neg(\sigma).\]
    \end{rem}
    We prove that the \emph{flag inversion} statistic is equidistributed with
    the \emph{flag-major index} over the \emph{colored permutation group} $G_{r,n}$:
    \begin{thm}\label{finvt}
        \[ \sum_{g \in G_{r,n}} q^{flag-major_F(g)}=
            \sum_{g \in G_{r,n}} q^{finv_F(g)}.\]
    \end{thm}
    \noindent{\bf Proof.} By using the proof of Lemma~\ref{tpc}, we can conclude:
    \begin{eqnarray*}
        \sum_{g \in G_{r,n}} q^{finv_F(g)}&=&
            \sum_{(z,id) \in Color^n_r}q^{csum((z,\pi))}\sum_{\pi \in S_n}
                q^{r \cdot inv(\pi)}\\
        &=& \sum_{(z,id) \in Color^n_r}q^{csum((z,id))}\sum_{\pi \in S_n}
                q^{r \cdot maj(\pi)}\\
        &=& \sum_{g \in G_{r,n}} q^{flag-major_F(\sigma)}
        \quad \bullet
    \end{eqnarray*}

    \begin{rem} There exists another direct proof of
    Theorem~\ref{finvt}. One can show by induction on $n$ that $a_n :=\sum\limits_{g \in
    G_{r,n}} q^{finv_F(g)}$ satisfies:
    \[a_n =a_{n-1} [r]_q[n]_{q^r}=a_n[rn]_q,\]
    and therefore
        \[\sum\limits_{g \in G_{r,n}} q^{finv_F(g)}=
            \prod_{i=1}^n[ri]_q.\]
        According to Theorem~\ref{rnmaj} we can conclude:
        \[ \sum_{g \in G_{r,n}} q^{flag-major(g)}=
            \sum_{g \in G_{r,n}} q^{finv_F(g)}.\]
    \end{rem}

    \begin{cor} The \emph{flag-inversion} statistic is
        equidistributed with the \emph{length} statistic over the
        \emph{signed permutation groups} $B_n$.
    \end{cor}
    \noindent{\bf Proof.} Combine to the Theorem~\ref{finvt} (in the
    case of $r=2$), with Proposition~\ref{fmajorder} and
    Theorem~\ref{fmajlen}.\\

    By using the \emph{flag-inversion} statistic we can generalize the Foata-Sch\"{u}tzenberger
    Theorem~\ref{FS1} to the \emph{colored permutation groups}:
    \begin{thm}\label{FSR}
        \[ \sum_{g \in G_{r,n}} q^{flag-major_F(g)}
            t^{ldes_F(\bar{g})}=\sum_{g \in G_{r,n}} q^{finv_F(g)}
            t^{ldes_F(\bar{g})},\]
        where $\bar{g}:=(z,\pi^{-1})$.
    \end{thm}
    \noindent{\bf Proof.} By using Lemma~\ref{tpc} and Theorem~\ref{FS2} we get:
     \begin{eqnarray*}
            \sum\limits_{g \in G_{r,n}}
                q^{flag-major_F(g)}t^{ldes_F(\bar{g})}
            &=& \sum\limits_{(z,\pi) \in G_{r,n}}
                    q^{flag-major_F((z,\pi))}t^{ldes_F((z,\pi^{-1}))}\\
            &=& \sum\limits_{(z,id) \in
                Color^n_r}q^{csum((z,id))}t^{csum((z,id))}
                \sum\limits_{\pi \in S_n} q^{r \cdot maj(\pi)}t^{des(\pi^{-1})}\\
            &=& \sum\limits_{(z,id) \in Color^n_r}
                q^{csum((z,id))}t^{csum((z,id))}
                    \sum\limits_{\pi \in S_n} q^{r \cdot inv(\pi)}t^{des(\pi^{-1})}\\
            &=&   \sum\limits_{(z,\pi) \in G_{r,n}}
                    q^{finv_F(g)}t^{ldes_F(\bar{g})} \quad \bullet
    \end{eqnarray*}
    In the case of $r=2$ by using \emph{flag-inversion} statistic, we can get the following
    equality, which is a type B generalization of the Foata-Sch\"{u}tzenberger
    Theorem~\ref{FS2} :
    \begin{thm}\label{GFS2}
    \[\sum_{\sigma \in B_n} q^{flag-major_F(\sigma^{-1})}t^{flag-major_F(\sigma)}=
        \sum_{\sigma \in B_n} q^{finv_F(\sigma)}t^{flag-major_F(\sigma)}.\]
    \end{thm}
    \noindent{\bf Proof.} Using Lemma~\ref{tpc}, Theorem~\ref{FS2}, and the fact that in the case
    of \\ $r=2, \ \tau=\tau^{-1} \in Color^n_2$ we get:
    \begin{eqnarray*}
            \sum\limits_{\sigma \in B_n}
                q^{flag-major_F(\sigma^{-1})}t^{flag-major_F(\sigma)}
            &=& \sum\limits_{\tau \in
                Color^n_2}q^{neg(\tau)}t^{neg(\tau)}
                \sum\limits_{\pi \in S_n} q^{2maj(\pi^{-1})}t^{2maj(\pi)}\\
            &=& \sum\limits_{\tau \in
                Color^n_2}q^{neg(\tau)}t^{neg(\tau)}
                \sum\limits_{\pi \in S_n} q^{2inv(\pi)}t^{2maj(\pi)}\\
            &=&   \sum\limits_{\sigma \in B_n}
                    q^{finv_F(\sigma)}t^{flag-major_F(\sigma)} \quad \bullet
    \end{eqnarray*}

    \medskip

    \subsection{Algebraic Interpretation}
        In this subsection we give an algebraic interpretation
        to the \emph{flag-inversion} statistic, as length in the
        group $S_{r\cdot n}$ with respect to a new set of generators.
        We open this subsection with the definition of the new
        generators:

        \begin{df} Let $r\geq 1,\ n \geq 1 $ be integers. Define the
        set of \emph{two dimensional Coxeter generators} of the \emph{Symmetric group} $S_{r\cdot n}$ by :
            \[C_2:=\{d_{i,j,R}  \ | 1\leq i \leq n,\ 0 \leq j \leq r-2 \}\bigcup
            \{d_{i,j,C} \ | 1\leq i \leq n-1,\ 0 \leq j \leq r-1 \},\]
            where:
                \[d_{i,j,R}:=[1,\ldots,1^{[r-1]},\ldots,i^{[j-1]},i^{[j+1]},i^{[j]},i^{[j+2]},\ldots,n,\ldots,
                n^{[r-1]}],\]
            and:
                \[d_{i,j,C}:=[\ldots,i^{[j-1]},(i+1)^{[j]},i^{[j+1]},\ldots,(i+1)^{[j-1]},i^{[j]},(i+1)^{[j+1]},\ldots].\]
        \end{df}

        There is an easy geometrical way to understand these
        generators; represent a permutation $\pi \in S_{r\cdot n}$ by a
        matrix $A=(a_{ij})_{i=1,\ j=0}^{i=n,\ j=r-1} \in M_{n \times r}$ , where $\forall i \in [1,n], k \in [1,r], \
        a_{i,k}=\pi(ir+k-1)$. For example, the identity permutation
        can be represented by the following matrix:
        \[ \left(
            \begin{array}{cccccc}
                1 &   1^{[1]}&    \ldots  &   1^{[r-2]} &   1^{[r-1]} \\
                2 &   2^{[1]}&    \ldots  &   2^{[r-2]} &   2^{[r-1]} \\
                \vdots &   \vdots   &   \ddots   &  \vdots   &   \vdots     \\
                n-1 &   (n-1)^{[1]}&    \ldots  &   (n-1)^{[r-2]} &   (n-1)^{[r-1]} \\
                 n &   n^{[1]}&    \ldots  &   n^{[r-2]} &   n^{[r-1]} \\
            \end{array} \right).\]
        The generator $d_{i,j,R}$ swaps two consecutive numbers in the same
        row, while $d_{i,j,C}$ swaps two consecutive numbers in the same column.
        \begin{eqnarray*}
          \left(
            \begin{array}{cccccc}
                1  \stackrel{\rm d_{1,0,R}}{\leftrightarrow} &
                1^{[1]}  \stackrel{\rm d_{1,1,R}}{\leftrightarrow} &
                \ldots   \stackrel{\rm d_{1,r-3,R}}{\leftrightarrow}&
                1^{[1]}  \stackrel{\rm d_{1,r-2,R}}{\leftrightarrow} &
                1^{[r-1]}\\ \\
                \updownarrow_{d_{1,0,C}}  &
                \updownarrow_{d_{1,1,C}}  & &
                 \updownarrow_{d_{1,r-2,C}}&
                \updownarrow_{d_{1,r-1,C}}&\\ \\
                2  \stackrel{\rm d_{2,0,R}}{\leftrightarrow} &
                2^{[1]}  \stackrel{\rm d_{2,1,R}}{\leftrightarrow} &
                \ldots   \stackrel{\rm d_{2,r-3,R}}{\leftrightarrow}&
                2^{[r-2]}  \stackrel{\rm d_{2,r-2,R}}{\leftrightarrow} &
                2^{[r-1]}  \\ \\

                \updownarrow_{\rm d_{2,0,C}}  &
                \updownarrow_{\rm d_{2,1,C}}  & &
                \updownarrow_{\rm d_{2,r-2,C}}  &
                \updownarrow_{\rm d_{2,r-1,C}} \\ \\

                \vdots &   \vdots   &   \ddots   &  \vdots &\vdots
                \\ \\
                \updownarrow_{\rm d_{n-1,0,C}}  &
                \updownarrow_{\rm d_{n-1,1,C}}  & &
                \updownarrow_{\rm d_{n-1,r-2,C}} &
                \updownarrow_{\rm d_{n-1,r-1,C}}& \\ \\
                n  \stackrel{\rm d_{n,0,R}}{\leftrightarrow} &
                n^{[1]}  \stackrel{\rm d_{n,1,R}}{\leftrightarrow} &
                \ldots   \stackrel{\rm d_{n,r-3,R}}{\leftrightarrow}&
                n^{[r-2]}  \stackrel{\rm d_{n,r-2,R}}{\leftrightarrow} &
                n^{[r-1]}  \\ \\
            \end{array} \right)
          \end{eqnarray*}
          \begin{center}
            \emph{figure 1. The Two dimensional Coxeter generators
            $d_{i,j,k}, \ k\in \{R,C\}$.}
          \end{center}
          \begin{rem} It is easy to check that for $r=1$
            the \emph{two dimensional Coxeter generators} $d_{i,0,C}$ are
            equal to the usual \emph{Coxeter generators} (Definition~\ref{CoxeterA}).
          \end{rem}

          We now define a bijection $\Gamma: G_{r,n} \to S_{r\cdot n}$, and calculate the length of $\Gamma(g)$
          according to the \emph{two dimensional Coxeter generators}:

          \begin{df} Define a bijection $\Gamma(g): G_{r,n} \to S_{r\cdot n}$ as follows: for $g \in G_{r,n}$
            \begin{eqnarray*}
                \Gamma(g)(i^{[k]}):=
                    \left\{%
                        \begin{array}{ll}
                            \pi(i)^{[z_i]}, & \hbox{$k=0$;} \\
                            \pi(i)^{[k-1]}, &  \hbox{$1 \leq k \leq z_i$;} \\
                            \pi(i)^{[k]}, & \hbox{$z_i < k \leq r-1$;} \\
                        \end{array}%
                    \right.
                , \ \forall k \in [0,r-1], i \in [1,n].
            \end{eqnarray*}
            Here we identify $[1,\ldots,r\cdot n]$ with
            $[1,\ldots,1^{[r-1]},\ldots,n,\ldots,n^{[r-1]}]$.
          \end{df}
    \begin{exa} Let $n=3,\ r=4,\ (z,\pi)=((0,3,2),[2,1,3])$. Then $\Gamma:
        G_{4,3} \to S_{12}$:
        \[\Gamma(g)=[2,2^{[1]},2^{[2]},2^{[3]},1^{[3]},1,1^{[1]},1^{[2]},3^{[2]},3,3^{[1]},3^{[3]}]=
            \left(
                \begin{array}{cccc}
                    2  &   2^{[1]}   &    2^{[2]}   &   2^{[3]} \\
                    1^{[3]}   &   1   &   1^{[1]}   &   1^{[2]} \\
                    3^{[2]}  &   3   &    3^{[1]}   &   3^{[3]}
                \end{array}
            \right).\]

    \end{exa}
    \begin{df} Let $g \in G_{r,n}$. Define the length of
        $g$ according to the \emph{two dimensional Coxeter
        generators} to be:
        \[\ell_D(g)=\min\{t \in \mathbf{N} \ | \
        \Gamma(g)=d_{i_1,j_1,k_1} \ldots d_{i_t,j_t,k_t},
        \mbox{ for some } i_s, \ j_s,\ k_s\}.\]
    \end{df}
    Our main result here is the next theorem.
    \begin{thm} Let $g \in G_{r,n}$. Then
        \[\ell_D(g)=r \cdot inv(g)+csum(g)=finv(g).\]
    \end{thm}

    \noindent{\bf Proof.} We give an algorithm which proves the
    theorem. Consider $\Gamma(g)$ as a matrix $A \in M_{n
    \times r}$. We do the following steps to create $A$ from the
    identity matrix using the generators $d_{i,j,k}$:
    \begin{enumerate}
        \item View every column of the matrix as a
            $\pi \in S_n$ use $d_{i,j,R}$ to generate it from the identity
            matrix.
            \[ \left(
                \begin{array}{ccc}
                    1   &   \ldots  &   1^{[r-1]}\\
                    \vdots  &   \ddots  &   \vdots \\
                    n   &   \ldots  &   n^{[r-1]}
                \end{array} \right)
                    \longrightarrow
              \left(
               \begin{array}{ccc}
                    \pi(1)   &   \ldots  &   \pi(1)^{[r-1]}\\
                    \vdots  &   \ddots  &   \vdots \\
                    \pi(n)   &   \ldots  &   \pi(n)^{[r-1]}
               \end{array} \right)\]
        \item In every row $j \in [1,n]$  move $j^{z_j}$ from
            column $z_j$  to column $1$ using $d_{i,j,C}$.
    \end{enumerate}
    \[
        \left(
        \begin{array}{ccc}
            \pi(1)   &   \ldots  &   \pi(1)^{[r-1]}\\
            \vdots  &   \ddots  &   \vdots \\
            \pi(n)   &   \ldots  &   \pi(n)^{[r-1]}
        \end{array} \right)
         \longrightarrow
        \left(
        \begin{array}{ccc}
            \pi(1)^{[z_1]}   &   \ldots  &   \pi(1)^{[r-1]}\\
            \vdots  &   \ddots  &   \vdots \\
            \pi(n)^{[z_n)}   &   \ldots  &   \pi(n)^{[r-1]}
        \end{array} \right)
        .\]
        The first step of the algorithm uses $r \cdot inv(\pi)$
        generators (to organize every column we use $inv(\pi)$  \emph{Coxeter
        generators} of $S_n$)
        and the second step uses $csum(g)$
        generators (in every column $i$ we move $\pi(i)^{[z_i]}$ from
        column $z_i$ to the first column, this action uses $z_i$
        generators). We conclude that $finv(g)\geq \ell_D(g)$. Now
        we need to prove that $\ell_D(g)\geq finv(g)$. We look at
        $\Gamma(g)$ as a matrix $(a_{ij})_{i=1,j=0}^{i=n,j=r-1}$, where
        $a_{ij}=\pi(i)^{[x_{ij}]}, \ x_{ij} \in [0,r-1]$. It easy to
        see that:
        \begin{eqnarray*}
            \ell_D(g)\geq
            r\cdot \ell(\pi)+\frac{1}{2}\sum_{i=1}^{n}\sum_{j=0}^{r-1}|x_{ij}-j|&=&
            r\cdot \ell(\pi)+\frac{1}{2}\sum_{i=1}^{n}( z_i+\sum_{j=0}^{z_i}1)\\
            &=&r\cdot \ell(\pi)+\sum_{i=1}^{n}z_i=finv(g).
        \end{eqnarray*}
        therefore;
        \[\ell_D(g)=finv(g) \ \ \bullet \]

\section{Acknowledgments}
I wish to thank my supervisors Prof. Ron M. Adin and Prof. Yuval
Roichman for their professional guidance, help and patience. I
have learned a lot from them.\\ \\
I am sincerely grateful to Taly Barak and her family for their love,
support and patience.\\ \\ I would like to thank all those who
helped and supported me during this work.


\begin{thebibliography}{99}

\bibitem{ABR}
R.\ M.\ Adin, F.\ Brenti and Y.\ Roichman, {\it Descent numbers
and major indices for the hyperoctahedral group}, Special issue in
honor of Dominique Foata's 65th birthday (Philadelphia, PA 2000),
Adv.\ Appl.\ Math.~{\bf 27} (2001), 210--224.

\bibitem{ARG}
R.\ M.\ Adin, I.\ M. Gessel and Y.\ Roichman {\em Signed
Mahonians}, J.\ Combin.\ Theory (ser.\ A)  ~{\bf 109} (2005), 25--43


\bibitem{AR}
R.\ M.\ Adin and Y.\ Roichman, {\em The flag major index and group
actions on polynomial rings}, Europ.\ J.\ Combin.~{\bf 22} (2001),
431--446.

\bibitem{BE}
E. \ Bagno, {\em Euler-Mahonian parameters on colored permutations groups}, S´minaire Lotharingien de Combinatoire ~{\bf 51} (2004), Article B51f .

\bibitem{BB}
A.\ Bj\"{o}rner and F.\ Brenti, {\em Combinatorics of Coxeter
Groups}, Graduate Texts in Mathematics, Springer-Verlag, to
appear.

\bibitem{Den}
M.\ Denert, {\em The genus zeta function of hereditary orders in
central simple algebras over global fields}, Math. Comp. {\bf 54}
(1990), 449-465.

\bibitem{DF}
J.\ D\'{e}sarmenien and D.\ Foata, {\em The signed Eulerian
numbers}, Discrete Math. {\bf 99} (1992),49 -58.

\bibitem{FS}
D.\ Foata and M.\ P.\ Sch\"utzenberger, {\em Major index and
inversion number of permutations}, Math.\ Nachr.~{\bf 83} (1978),
143--159.

\bibitem{FZ}
D.\ Foata and D.\ Zeilberger, {\em Denert's permutation statistic
is indeed Euler-Mahonian}, Studies in Appl. Math. {\bf 83} (1990),
31-59.

\bibitem{G}
I. \ Gessel, {\it Generating Functions and Enumeration of
Sequences}, Ph. D. Thesis, MIT, (1977).

\bibitem{Hu}
J.\ E.\ Humphreys, {\em Reflection Groups and Coxeter Groups},
Cambridge Studies in Advanced Mathematics, no.~29, Cambridge
Univ.\ Press, Cambridge, 1990.

\bibitem{MM}
P.\ A.\ MacMahon, {\em Combinatory Analysis I-II}, Cambridge
Univ.\ Press, London/New-York, 1916. (Reprinted by Chelsea,
New-York, 1960.)

\bibitem{ST}
E.\ Steingr\'imsson, {\em Permutation statistics of indexed permutations}, Europ.
J. Combin. {\bf 15} (1994), 187-–205.

\bibitem{Wachs}
M.\ Wachs, {\em An involution for signed Eulerian numbers}, Discrete
Math.\ {\bf 99} (1992), 59--62.



\end{thebibliography}
\end{document}